\documentclass[11pt]{amsart}
\usepackage[a4paper,top=4cm,bottom=3cm,left=3cm,right=3cm]{geometry}
\usepackage[T1]{fontenc}
\usepackage{bm}
\usepackage{amsmath,amssymb,amsfonts,amsthm} 
\usepackage{hyperref}
\usepackage{graphicx}
\usepackage{dirtytalk}
\usepackage{pgfplots}
\usepackage{yfonts}
\usepackage{url}
\usepackage{enumitem}
\usepackage{esint}
\usepackage{multirow}
\usepackage{blindtext}
\newcommand{\numberset}{\mathbb}

\newcommand{\R}{\numberset{R}}
\newcommand{\M}{\textup M}
\newcommand{\g}{\textup g}

\newcommand{\cs}{\textup{c}_\textup{s}}
\newcommand{\cper}{\textup{c}_\textup{p}}

\newcommand{\cin}{\textup{c}_{\textup{in}}}

\newtheorem{thm}{Theorem}[section]
\newtheorem{rem}[thm]{Remark}
\newtheorem{lem}[thm]{Lemma}
\newtheorem{dfn}[thm]{Definition}
\newtheorem{cor}[thm]{Corollary}
\newtheorem{prp}[thm]{Proposition}

\title[Volume preserving mean curvature flow]{Volume preserving mean curvature flow of round surfaces in asymptotically flat spaces}

\author{Carlo Sinestrari}
\address{(Carlo Sinestrari) Dipartimento di Matematica, Universit\`a degli Studi di Roma "Tor Vergata", Via della Ricerca Scientifica, 00133, Roma, Italy}
\author{Jacopo Tenan}
\address{(Jacopo Tenan) Dipartimento di Matematica, Universit\`a degli Studi di Roma "Tor Vergata", Via della Ricerca Scientifica, 00133, Roma, Italy}

\date{}
\numberwithin{equation}{section}
\pgfplotsset{compat=1.18}
\begin{document}

\begin{abstract}
We study the volume preserving mean curvature flow of a surface immersed in an asymptotically flat $3$-manifold modeling an isolated gravitating system in General Relativity. We show that, if the ambient manifold has positive ADM mass and the initial surface is round in a suitable sense, then the flow exists for all times and converges smoothly to a stable CMC surface. This extends to the asymptotically flat setting a classical result by Huisken-Yau (Invent. Math. 1996) and allows to construct a CMC foliation of the outer part of the manifold by an alternative approach to the ones 
by Nerz (Calc. Var. PDE, 2015) or by Eichmair-Koerber (J. Diff. Geometry, 2024).
\end{abstract}

\maketitle

\section{Introduction}

In this paper we study the evolution by mean curvature of closed surfaces in a smooth Riemannian manifold which is asymptotically flat, according to the following definition.

\begin{dfn}\label{asymflat}A  complete Riemannian 3-manifold $(\M,\overline\g)$ is called $C_{\frac12+\delta}^2$-\emph{asymptotically flat} for some $\delta\in(0,\frac12]$ if there exists a compact subset $\textup C\subset \M$, a constant $\overline c>0$ and a diffeomorphism $\vec x:\M\setminus \textup C\to \R^3\setminus \overline{\mathbb{B}}_1(0)$ such that 
\begin{equation}\label{asymptc8}
|\overline \g_{\alpha\beta}-\overline \g^e_{\alpha\beta}|+|\vec x|\left| \partial_\gamma \, \overline \g_{\alpha\beta}\right|+|\vec x|^2\left|\partial_\gamma\partial_\eta \, \overline \g_{\alpha\beta}\right|\leq \overline c|\vec x|^{-\frac12-\delta}
\end{equation}
and in addition the scalar curvature $\overline{\textup{S}}$ satisfies
$|\overline{\textup{S}}|\leq \overline c|\vec x|^{-3-\delta}.$ Here $\overline \g^e$ denotes the Euclidean metric, and
the partial derivatives 
and the norms are the Euclidean ones.
\end{dfn}

Manifolds of this kind have been extensively studied in General Relativity, since they occur as spacelike time slices of Lorentzian manifolds modeling isolated gravitating systems. A typical example is the well-known Schwarzschild metric of mass $m$, given by
\begin{equation}
\bar g^S_{\alpha \beta} = \left( 1+\frac{m}{2|\vec x|} \right)^4 \overline \g^e_{\alpha\beta},
\end{equation}
which fulfills the above definition with $\delta=\frac 12$.

In this paper we study the \textit{volume preserving mean curvature flow} (VPMCF) of surfaces immersed in $\M$, namely we consider time dependent immersions $F:\Sigma \times [0,T) \to \M$, with $\Sigma$ a closed $2$-surface, which evolve according to
\begin{equation}\label{vpmcf}
\frac{\partial F}{\partial t}(p,t) = -[H(p,t)-h(t)] \,\nu (p,t).
\end{equation}
Here $H$ and $\nu$ are the mean curvature and the unit normal, while $h(t)$ is the average of the mean curvature on $\Sigma$ at time $t$. Compared with the standard mean curvature flow, where the $h(t)$ term is absent, this evolution is such that the volume of the region enclosed by the surface remains constant in time. A remarkable feature of the flow is that the area of the surface decreases with time, with strict monotonicity unless the curvature is constant.
Convergence results of the flow to a CMC profile for suitable classes of initial data, both in Euclidean and in Riemannian ambient spaces, have been obtained by many authors through the years, we recall for instance \cite{alfr,cabezasrivas,essi,huisken1987,jmos,li}.

In the context of mathematical relativity, VPMCF was first studied by Huisken and Yau \cite{huiskenyau} who proved the following result.
\begin{thm}
Let $(\M,\overline\g)$ a $C^4_2$-\emph{asymptotically Schwarzschild} $3$-manifold, in the sense that \eqref{asymptc8} is replaced by
\begin{equation}\label{asymptS}
|\overline \g-\overline \g^S|+\sum_{k=1}^4|\vec x|^k\left| \partial^{|k|} \, (\overline \g_{\alpha\beta}-\overline \g^S_{\alpha\beta})\right| \leq \overline c|\vec x|^{-2}
\end{equation}
where $\overline \g^S$ is the Schwarzschild metric for some $m>0$ and  $\partial^{|k|}$ denote derivatives of order $k$. 
Let $\Sigma_t=F(\Sigma,t)$ be the solution of the flow \eqref{vpmcf} with initial data given by the Euclidean coordinate sphere $\mathbb{S}_r(0)$, for a large enough radius $r>0$. Then $\Sigma_t$ exists for all $t \in [0,\infty)$ and converges smoothly to a stable CMC-surface $\Sigma^r_\infty \subset M$ as $t \to +\infty$.
\end{thm}

Huisken and Yau then showed that the union of the surfaces $\Sigma^r_\infty$ as $r \in (r_0,+\infty)$ forms a foliation by stable CMC-surfaces of the outer part of $\M$, which is uniquely determined. Such a foliation is of interest for the physical model because it defines a canonical system of radial coordinates and allows for a geometric definition of the center of mass. During the years, various authors have given alternative constructions of the CMC foliation and have weakened the assumptions, as in the papers by Ye \cite{ye}, Metzger \cite{metzger} and Huang \cite{huang1}. An important reference for our purposes is the work of Nerz \cite{nerz1}, who first proved the result under the $C_{\frac12+\delta}^2$-flatness decay assumption, which he then showed to be optimal in \cite{nerzJFA}. We also mention the recent paper by Eichmair-Koerber \cite{eichmairkoerber}, who have constructed the CMC foliation in general dimensions and have given a survey of the previous results.

In the above mentioned papers that followed \cite{huiskenyau}, the existence of the foliation was obtained by more classical tools of elliptic theory instead of using curvature flow evolutions. This has motivated our present work, aimed at extending the analysis of the VPMCF by Huisken and Yau to general  $C_{\frac12+\delta}^2$-asymptotically flat manifolds with positive ADM-mass. Our main result, Corollary \ref{mainthPart1}, states that the flow exists for all times and converges to a CMC limit provided the initial data is enough round in a suitable sense. In particular, if we add a mild symmetry assumption on the ambient metric, we obtain the following statement about the evolution of the coordinate spheres.

\begin{thm}\label{mainspheres}
Let $(\M,\overline\g)$ be a $C_{\frac12+\delta}^2$-asymptotically flat 3-manifold with $m_\textup{ADM}>0$, satisfying the $C_{1+\delta}^1$-Regge-Teitelboim conditions in Definition \ref{weakrtcond}. Then the solution of the flow \eqref{vpmcf} starting from a Euclidean coordinate sphere with large enough radius exists for all times and converges to a stable CMC-surface as $t \to +\infty$.
\end{thm}

The surfaces obtained as limits of the flow coincide with the ones constructed in \cite{nerz1} and therefore our analysis provides an alternative approach to the existence of the CMC-foliation in the $C_{\frac12+\delta}^2$-{asymptotically flat} setting, giving also
a quantitative dynamical description of the stability property of the surfaces. On the other hand, we remark that our result, in contrast with \cite{nerz1}, requires the hypothesis of positive mass. In fact, if the mass is negative,  the surfaces of the foliation are no longer area minimizing under volume constraint, and one expects  that neighbouring surfaces generically drift away with the flow.

As in \cite{huiskenyau}, the main idea in our construction is to define a suitable class of immersed surfaces, called {\em round surfaces}, which are close to Euclidean spheres and are well-centered with respect to the Euclidean coordinates. The properties of a round surface are quantitatively expressed in terms of an approximate radius $\sigma$ and become more restrictive as $\sigma \to +\infty$ according to the increasingly Euclidean behaviour of the ambient metric. The central step in the proof is to show that, if the parameters of the class are suitably chosen and the radius is sufficiently large, the solution of the flow remains round as long at it exists. The long time existence and the convergence to a CMC limit then follow from standard arguments.

To prove that the solution of the flow stays inside the roundness class, we need a different strategy from the one in \cite{huiskenyau},
because of the weaker assumptions on the ambient metric which prevent the use of the maximum principle to obtain curvature estimates.
Following an idea of Metzger \cite{metzger}, we study the invariance of integral norms of the curvature, instead of pointwise bounds as in \cite{huiskenyau}.
Conversely, the integral estimates allow to obtain a pointwise control on the curvature by using some regularity results from the literature: in particular 
the rigidity estimate for nearly umbilical surfaces by De Lellis and M\"uller \cite{delellismuller1}, together with its $L^p$ version by Perez \cite{perez} and the bootstrap for the second fundamental form of nearly CMC-surfaces by Nerz \cite{nerz1}. We observe that integral estimates have been often used to obtain convergence of the mean curvature flow, starting from the pioneering paper \cite{huisken1}. Here, however, we do not need to employ directly a Stampacchia iteration technique to obtain $L^\infty$ estimates, since this step is implicitly contained in the result by Nerz, who uses a (simpler) integral iteration procedure in the elliptic setting.
 
As in previous works, an important part in our procedure is played by the spectral properties of the stability operator associated to the evolving surface, which we denote by $L$. In particular, the positivity of the smallest eigenvalue of $L$, restricted on the functions with zero mean value, is related to the decay of the $L^2$-norm of $H-h$ along the flow, which in turn controls the possible drift of the barycenter of the surface. The spectral properties of $L$ have been analyzed in detail  in \cite{cederbaum,nerz1} in the case of round surfaces which are CMC, showing in particular that the smallest eigenvalue has the same sign as the ADM-mass of $\M$. Our evolving surfaces, however, are not CMC but only satisfy a smallness assumption on the oscillation of $H$. This gives rise to some non negligible remainder terms and the positivity of $L$ may fail. On the other hand, for our purposes, we only need to apply $L$ to the deformation $H-h$ induced by the flow, and we prove that in this case the remainder terms admit a better estimate thanks to symmetry properties. This leads to an exponential decay in $L^2$-norm for the speed of the flow, see Proposition \ref{evolutionH-h}, which gives the desired control on the barycenter if the so-called translational part of $H$, see Definition \ref{spectral}, is enough small at the initial time. It is then easy to check that the Euclidean coordinate spheres with large enough radius satisfy this smallness assumption if the ambient space satisfies a weak Regge-Teitelboim condition, and this yields Theorem \ref{mainspheres}.

We conclude by describing further  related works in the literature. Other convergence results for constrained mean curvature flows in asymptotically Schwarzschild manifolds were obtained in \cite{corvinowu} and recently in \cite{guilisun}. An interesting modification of the approach by Huisken-Yau has been given by Cederbaum-Sakovich \cite{cederbaum} by considering the so-called space-time mean curvature (STMC): by using a strategy similar to \cite{nerz1}, the authors construct a constant STMC foliation which allows to treat certain cases, described in \cite{cederbaumnerz}, where the center of mass via the CMC-foliation is not well-defined. In this context, the second author \cite{tenan} has considered the volume preserving STMC-flow and has shown that it drives the CMC leaves constructed in the present paper to a constant STMC limit, providing an alternative construction of the foliation in the positive energy case. Very recently, Kr\"oncke-Wolff \cite{krwo} have used an area preserving null mean curvature flow to construct foliations of asymptotically Schwarzschildean lightcones by surfaces of constant STMC. \\
 
\noindent\textbf{Acknowledgments:} Both authors have been supported by MUR Excellence Department Project MatMod@TOV  CUP E83C23000330006 and by the MUR Prin 2022 Project  ``Contemporary perspectives on geometry and gravity'' CUP E53D23005750006. C.S.
is a member of the group GNAMPA of INdAM (Istituto Nazionale di Alta Matematica). C.S. is grateful to Gerhard Huisken for inspiring discussions and suggestions.
\bigskip

\section{Preliminaries}

\subsection{Definitions and basic properties}

Throughout the paper, $(\M,\overline\g)$ will be an asymptotically flat $3$-manifold as in Definition \ref{asymflat}. 
We denote respectively by $\overline \Gamma^\alpha_{\beta \gamma}$ the Christoffel symbols, by $\overline \nabla$ the Riemannian connection, by $\overline {\textup{Rm}}, \overline{\textup{Ric}}$ the Riemann and Ricci curvature tensors and by $\overline{\textup{S}}$ the scalar curvature on $(\M,\overline\g)$.
We call \emph{Euclidean coordinate spheres} the surfaces of the form $\vec x^{-1}(\mathbb{S}_R(0))$ for some $R>0$; by an abuse of notation,
we denote them simply by $\mathbb S_R(0)$. 

The \emph{ADM-mass} of $(\M,\overline\g)$, first introduced in \cite{adm}, is defined as 
\begin{equation}\label{admenergy}
\overline m_{\textup{ADM}}:=\lim_{R\to\infty}\frac{1}{16\pi}\int_{\mathbb S_R( 0)}
\left(  \partial_\alpha  \overline \g_{\alpha\beta} -  \partial_\beta  \overline \g_{\alpha\alpha} \right) \nu^\beta \ d\mu,
\end{equation}
where $\nu$ is the unit normal vector and $d\mu$ is the measure on $\mathbb S_R( 0)$ induced by $\overline \g$.

We also introduce a weakened form of a symmetry assumption originally stated by Regge and Teitelboim \cite{reggeteitelboim}. We remark that this property will only be used in Theorem \ref{mainspheres} and Lemma \ref{rtspheres}, while it is not needed in the other results of the paper.

\begin{dfn}[$C_{1+\delta}^1$-Regge-Teitelboim conditions]\label{weakrtcond} Let $(\M,\overline \g)$ be a $C_{\frac12+\delta}^2$-asymptotically flat manifold. We say that this manifold satisfies the $C_{1+\delta}^2$-Regge-Teitelboim conditions if there exists $\overline c>0$ such that 
\begin{equation}\label{wrtcond363}
\left|\overline \g_{\alpha \beta}(\overline x)-\overline\g_{\alpha \beta}(-\overline x)\right|+|\overline x|\left|\overline\Gamma_{\alpha \beta}^\gamma(\overline x)+\overline\Gamma_{\alpha \beta}^\gamma(-\overline x)\right|\leq \frac{\overline c}{|\overline x|^{1+\delta}}
\end{equation}
for every $\overline x\in \M\setminus \textup C$.
\end{dfn}

In the following, we will consider smooth immersed $2$-surfaces $\iota:\Sigma\hookrightarrow \M\setminus \textup C$. We will always assume that $\Sigma$ is closed, connected and orientable. On $\Sigma$ we can consider the physical metric $g:=\iota^*\overline\g$ and the Euclidean metric $g^e:=\iota^*\overline \g^e$, where $\overline \g^e$ is induced by the diffeomorphism $\vec x$.  We will use the apex $e$ on each geometric quantity when it is computed with respect to the Euclidean metric, and we write $\Sigma^e$ as a short notation for $(\Sigma,g^e)$.
We use the latin indices $i,j,k,\dots \in \{1,2\}$ to denote the coordinates on $\Sigma$, while the coordinates in the ambient space are indicated with the greek letters $\alpha,\beta,\gamma,\dots \in \{1,2,3\}$.

We denote by $\nabla$ the connection on $\Sigma$, by $S$ the scalar curvature and by $d\mu$ the volume form. We denote by $\nu$ the unit outer normal, by $A=\{h_{ij}\}$ and by $H$, respectively, the second fundamental form and the mean curvature, and by $\kappa_i$, $i\in\{1,2\}$ the principal curvatures. In addition, we denote by  $\overset{\circ}{A}:= A - \frac{H}{2}g$ the traceless part of the second fundamental form. We recall the identities
$$
|A|^2=\kappa_1^2+\kappa_2^2, \qquad H^2=(\kappa_1+\kappa_2)^2, \qquad |\overset{\circ}{A}|^2=|A|^2-\frac 12 H^2= \frac 12 (\kappa_1-\kappa_2)^2.
$$
The property of asymptotic flatness allows to estimate the difference between these quantities and their counterparts computed with respect to the Euclidean metric. We collect the relevant properties in the following lemma, which can be proved by standard computations, see e.g. \cite[Lemma 1.5]{metzger}, \cite[Lemma 11]{cederbaum}.

\begin{lem}\label{metzger}
Let $\iota:\Sigma\hookrightarrow \M$ be a surface immersed in a $C_{\frac12+\delta}^2$-asymptotically flat 3-manifold $\M$. Then there exists  $C>0$, only depending on the constant $\bar c$ in \eqref{asymptc8}, such that
\begin{equation}
|g-g^e|_g \leq C |\vec x|^{-\frac 12-\delta}, \qquad |\Gamma^k_{ij}-(\Gamma^e)^k_{ij}| \leq C |\vec x|^{-\frac 32-\delta}
\end{equation}
\begin{equation}\label{dmu}
|d\mu-d\mu^e| \leq C |\vec x|^{-\frac 12-\delta}d\mu,
\end{equation}
\begin{equation}
|\nu-\nu^e|_g \leq C |\vec x|^{-\frac 12-\delta}, \qquad |\nabla \nu-\nabla^e \nu^e|_g \leq C |\vec x|^{-\frac 32-\delta},
\end{equation}
\begin{equation}
|A-A^e|_g \leq C \left(|\vec x|^{-\frac 32-\delta} +|\vec x|^{-\frac 12-\delta} |A^e| \right)
\end{equation}
\begin{equation}
|\nabla A-\nabla^e A^e|_g \leq C \left(|\vec x|^{-\frac 52-\delta} +|\vec x|^{-\frac 12-\delta} |\nabla^e A^e| \right).
\end{equation}
In addition, if $|A| \leq 10 |\vec x|^{-1}$, we have
\begin{equation}
|H-H^e| \leq C |\vec x|^{-\frac 32-\delta}, \qquad | \overset{\circ}{A} - \overset{\circ}{A^e}|_g \leq C |\vec x|^{-\frac 32-\delta}.
\end{equation}
\end{lem}

We introduce some more notation. We denote the average of the mean curvature by
\begin{equation*}
h:=\frac1{|\Sigma|}\int_\Sigma H\ d\mu,
\end{equation*}
where $|\Sigma|=\int_\Sigma d\mu$. In addition, we define the barycenter of $\Sigma$ as  
\begin{equation*}
\vec z_\Sigma:=\frac{1}{|\Sigma|}\int_\Sigma \iota \ d\mu,
\end{equation*}
that is, the average of the (Euclidean) position vector with respect to the physical metric. This definition differs slightly from the one in \cite{cederbaum,nerz1}, where the average is taken in the Euclidean metric; however, the two definitions have the same qualitative properties.

\begin{dfn}\label{defradii} Let $\iota:\Sigma\hookrightarrow \M$ be a surface immersed in an asymptotically flat 3-manifold $\M$. Then we set 
\begin{equation*}
r_\Sigma:=\min_{x\in\Sigma}|\vec x(\iota(x))|, \qquad R_\Sigma:=\max_{x\in\Sigma}|\vec x(\iota(x))|, \qquad \sigma_\Sigma:=\sqrt{\frac{|\Sigma|}{4\pi}}.
\end{equation*}
We call these values \emph{Euclidean radius, Euclidean diameter} and \emph{area radius}, respectively. 
\end{dfn}

As in \cite{cederbaum,nerz1} we consider the Sobolev norms on $\Sigma$ with a radius-dependent weight as follows
$$
\|f\|_{W^{0,p}(\Sigma)}:=\|f\|_{L^p(\Sigma)}, \qquad \|f\|_{W^{k+1,p}(\Sigma)}:=\|f\|_{L^p(\Sigma)}+\sigma_\Sigma \|\nabla f \|_{W^{k,p}(\Sigma)},
$$
for $p \in [1,\infty]$ and $k \geq 0$ integer. As usual, we use the notation $H^k=W^{k,2}$.

\subsection{Round surfaces} 
We now introduce a class of surfaces, which are close to a Euclidean sphere of radius $\sigma$ in a quantitative way measured by some parameters. The aim is to find a class which is invariant under the volume preserving mean curvature flow for an appropriate choice of the parameters and for large enough radius, or at least such that the possible deterioration of the parameters along the flow can be estimated. Other classes of round surfaces, which are related to the methods used there, have been introduced in \cite{huiskenyau,metzger,nerz1}

\begin{dfn}[Round surfaces]\label{roundsurf251} Let $(\M,\overline \g)$ be a $C_{\frac12+\delta}^2$-asymptotically flat manifold and let $\iota:\Sigma\hookrightarrow \M$ be a surface. \\
\textup{(i)} For a given approximate radius $\sigma>1$, and parameters $\eta,B_1,B_2>0$, we say that $\iota(\Sigma)$ is a \emph{round surface} in $(\M,\overline \g)$, and we write $\iota(\Sigma) \in
\mathcal{W}^\eta_\sigma(B_1,B_2)$,
if the following inequalities are satisfied:
\begin{equation}\label{controlon|A|}
|A(t)| < \sqrt{\frac 52}\sigma^{-1}, \qquad \kappa_i(t) >\frac{1}{2\sigma},
\end{equation}
\begin{equation}\label{radiicondition38}
\frac 72\pi\sigma^2< |\Sigma| \leq 5\pi\sigma^2, \qquad \frac34<\frac{r_\Sigma}{\sigma}< \frac{R_\Sigma}{\sigma}\leq \frac54,
\end{equation}
\begin{equation}\label{cond2defround}
\|\overset{\circ}{A}\|_{L^4(\Sigma)}<B_1\sigma^{-1-\delta}, 
\end{equation}
\begin{equation}\label{cond3defround}
\eta\sigma^{-4}\|H-h\|_{L^4(\Sigma)}^4+\|\nabla H\|_{L^4(\Sigma)}^4<B_2\sigma^{-8-4\delta}.
\end{equation}
\textup{(ii)}  For given $\sigma>1$ and  $\eta,B_1,B_2,B_\textup{cen}>0$, we say that $\iota(\Sigma)$ is a \emph{well-centered round surface}, and we write $\iota(\Sigma) \in
 \mathcal{B}^\eta_\sigma(B_1,B_2,B_\textup{cen})$ if it satisfies the above properties and in addition
\begin{equation}\label{cond2defroundprima}
|\vec z_\Sigma|<B_\textup{cen}\sigma^{1-\delta}.
\end{equation}
\end{dfn}
\medskip

At a first sight it may look redundant to introduce a further scale parameter $\sigma$, in addition to the ones in Definition \ref{defradii}, since condition \eqref{radiicondition38} implies that each of these values controls the others. However, this choice is practical in the study of the curvature flow evolution, where $\Sigma$ depends on time. Then the radii of Definition \ref{defradii} change with time, and it is convenient to describe the size of the evolving surface in terms of a fixed parameter $\sigma$.

\begin{rem} 
The decay rates in conditions \eqref{cond2defround}-\eqref{cond3defround} are modeled on the ones of the Euclidean coordinate spheres. In fact, by Lemma \ref{metzger}, it is easy to check that if $B_1,B_2$ are large enough, depending on $\bar c$ in \eqref{asymptc8}, then $\mathbb{S}_r(0)$ belongs to $\mathcal{B}^\eta_\sigma(B_1,B_2,B_\textup{cen})$ for $r$ large enough and $r/\sigma$ enough close to $1$. Conversely, we will see in Lemma \ref{cor1}(iv) that a round surface is close to a sphere in Euclidean coordinates; in particular, it is embedded and has genus zero. It is also easy to see that a round surface in our sense also belongs to the class of asymptotically concentric surfaces defined in \cite[Def. 4.3]{nerz1}. 
\end{rem}

\noindent
{\bf On the notation for the constants.} Throughout the paper, when deriving estimates on geometric quantities on a surface $\Sigma$, 
we denote by $C,C_1,C_2,\dots$ constants which only depend on properties of the ambient manifold, such as $\bar c, \delta$ in \eqref{asymptc8} or the mass $\overline m_{\textup{ADM}}$ and by $c,c_1,c_2,\dots$ constants which in addition depend on the constants $B_1,B_2,B_\textup{cen}$ in the previous definition. We say that a constant is universal if it is independent on any other parameter of our problem. As usual, the letters $c$ or $C$ will often denote constants which may change from one line to the other, but each time depending on the same parameters.
\medskip

Let us first observe some easy consequences of the properties of a round surface.

\begin{rem}\label{firstround}
{\rm (i)} Property \eqref{radiicondition38} implies that the three radii of Definition \ref{defradii} are comparable among each other and with $\sigma$.
In particular this property implies, because of the asymptotic flatness of $\M$ in \eqref{asymptc8}, the bound on the Riemann tensor
\begin{equation}\label{boundonriem}
 |\overline{\textup{Rm}}| \leq C(\bar c) \sigma^{-\frac 52 -\delta} \mbox{ on } \Sigma.
\end{equation}
{\rm (ii)} The Michael-Simon inequality in Euclidean space, together with the curvature bound in \eqref{controlon|A|} and the error estimates in Lemma \ref{metzger}, implies the existence of a universal Sobolev constant $c_S>0$ such that
\begin{equation}\label{sobolevineq}
\|\psi\|_{L^2(\Sigma)}\leq\frac{\cs}{\sigma}\|\psi\|_{W^{1,1}(\Sigma)}, \qquad \forall \ \psi\in W^{1,1}(\Sigma),
\end{equation}
on any round surface $\Sigma$ with radius $\sigma \geq \sigma_0=\sigma_0(\overline c,\delta)>0$. From this, the other Sobolev inequalities can be deduced. In particular (see e.g. Lemma 12 in \cite{cederbaum} and the references therein) we have, for every $p>2$,
\begin{equation}\label{othersob11}
\|\psi\|_{L^\infty(\Sigma)}\leq 2^{\frac{2(p-1)}{p-2}}\cs\sigma^{-\frac{2}{p}}\|\psi\|_{W^{1,p}(\Sigma)}, \qquad \forall \psi\in W^{1,p}(\Sigma),
\end{equation}
and also
\begin{equation}\label{sobH2}
\|\psi\|_{L^\infty(\Sigma)}\leq 32 \cs^2\sigma^{-1}\|\psi\|_{H^2(\Sigma)}, \qquad \forall \psi\in H^2(\Sigma).
\end{equation}
\end{rem}

In the next Lemma, similar to \cite[Proposition 4.4]{nerz1} and \cite[Proposition 1]{cederbaum}, we collect some properties of round surfaces which follow from the integral bounds \eqref{cond2defround} and \eqref{cond3defround} by applying the regularity results by De Lellis and M\"uller  \cite{delellismuller1} and by Nerz \cite{nerz1}.

\begin{lem}\label{cor1}
Let $(\M,\overline \g)$ be a $C_{\frac12+\delta}^2$-asymptotically flat manifold. Fix any $\eta,B_1,B_2>0$. Then there exists $\sigma_0=\sigma_0(B_1,B_2,\eta,\overline c,\delta)>0$ such that any surface $(\Sigma,g)$ which belongs to $\mathcal{W}_\sigma^\eta(B_1,B_2)$ for some $\sigma>\sigma_0$ satisfies the following properties.
\begin{enumerate}[label=\textup{(\roman*)}]
\item There exists a constant $c=c(B_2,\eta)>0$ such that  
\begin{equation}\label{LinfcontrolH}
\|H-h\|_{L^\infty}\leq c \sigma^{-\frac32-\delta}.
\end{equation}
\item There exists $c=c(B_1,B_2,\eta,\overline c)>0$ such that 
\begin{equation}\label{eqbis529}
\left|h-\frac2{\sigma_\Sigma}\right|\leq c \sigma^{-\frac32-\delta}.
\end{equation}
\item There exists a constant $B_\infty=B_\infty(B_1,B_2,\eta,\overline c,\delta)$ such that $\|\overset{\circ}{A}\|_{L^\infty(\Sigma)}\leq B_\infty\sigma^{-\frac32-\delta}$.
\item There exist constants $c=c(B_1,B_2,\eta,\overline c,\delta)$, $c_0=c_0(B_1,\overline c,\delta)$, $\vec z_0\in\mathbb R^3$, and a function $f:\mathbb{S}_{\sigma_\Sigma}(\vec z_0)\to \mathbb{R}$ such that 
\begin{equation}\label{eq2110}
\Sigma^e=\textup{graph}(f), \qquad \|f\|_{W^{2,\infty}}\leq c\sigma^{\frac12-\delta}, \qquad |\vec z_0-\vec z_\Sigma|\leq c_0\sigma^{\frac12-\delta}.
\end{equation}
\item There exists a constant $\cper=\cper(\overline c,\delta)$ such that 
\begin{equation*}
\|H-h\|_{L^4(\Sigma)}\leq \cper\|\overset{\circ}{A}\|_{L^4(\Sigma)}+\cper\sigma^{-1-\delta}.
\end{equation*}
\end{enumerate}
\end{lem}
We point out that the interest of part (v) comes from the fact that $\cper$ is independent on the constant $B_2$ which appears in \eqref{cond3defround}.

\begin{proof}
(i) The estimate follows from property  \eqref{cond3defround} and the Sobolev immersion \eqref{othersob11}. \medskip\\
(ii) We first observe that  \eqref{radiicondition38}, \eqref{cond2defround} and \eqref{cond3defround} imply,
\begin{equation}\label{L2norms}
\|\overset{\circ}{A}\|_{L^2(\Sigma)}<c_1\sigma^{-\frac12-\delta}, \qquad \| H-h\|_{H^1(\Sigma)}<c_2\sigma^{-\frac12-\delta},
\end{equation}
where $c_1=c_1(B_1)$ and $c_2=c_2(\eta,B_2)$. Then we recall the estimate by De Lellis and M\"uller in \cite{delellismuller1} for surfaces in Euclidean space, which states that
$$
\left\|A^e-\frac{1}{\sigma_{\Sigma^e}} g^e \right\|_{L^2(\Sigma^e)} \leq c_{DM} \left\|\overset{\circ}{A^e}\right\|_{L^2(\Sigma^e)},
$$ 
where $c_{DM}>0$ is a universal constant and $\sigma_{\Sigma^e}$ is defined by $4 \pi \sigma_{\Sigma^e}^2=|\Sigma^e|$. Passing to the physical metric, using Proposition \ref{metzger} and properties 
\eqref{controlon|A|}, 
\eqref{radiicondition38} we deduce
\begin{equation}\label{DMriem}
\left\|A-\frac{1}{\sigma_\Sigma} g \right\|_{L^2(\Sigma)} \leq c_{DM} \left\|\overset{\circ}{A}\right\|_{L^2(\Sigma)}+ C(\overline c) \sigma^{-\frac 12-\delta} \leq c \sigma^{-\frac 12-\delta},
\end{equation}
for some $c=c(B_1,\overline c)$. On the other hand, we have 
\begin{eqnarray*}
\sqrt{2\pi} \sigma_\Sigma\left|h-\frac2{\sigma_\Sigma}\right| & = &
 \left\|\frac{h}2g-\frac{1}{\sigma_\Sigma}g\right\|_{L^2(\Sigma)} \\
& \leq & \left\|\frac{h}2g-\frac{H}{2}g\right\|_{L^2(\Sigma)}+\left\|\frac{H}{2}g-A\right\|_{L^2(\Sigma)} +
\left\|A-\frac{1}{\sigma_\Sigma} g \right\|_{L^2(\Sigma)} \\
&= & \frac{1}{\sqrt 2} \left\|h-H \right\|_{L^2(\Sigma)}+  \left\|\overset{\circ}{A}\right\|_{L^2(\Sigma)} +
\left\|A-\frac{1}{\sigma_\Sigma} g \right\|_{L^2(\Sigma)}.
\end{eqnarray*}
Taking into account \eqref{radiicondition38},\eqref{L2norms} and \eqref{DMriem}, we obtain \eqref{eqbis529}. \medskip \\
(iii)
Assumptions \eqref{cond2defround} and \eqref{cond3defround}, together with the properties in Remark \ref{firstround}, allow to apply the bootstrap regularity for the second fundamental form by Nerz \cite[Proposition 4.1]{nerz1}, which give the assertion. \medskip \\
 (iv) The result of De Lellis-M\"uller quoted above \cite[Thm. 1.1]{delellismuller1} also gives the existence of a conformal parametrization 
$\Psi:\mathbb{S}_{\sigma_{\Sigma^e}}(\vec z_0)\to \Sigma^e$ for a suitable center $\vec z_0 \in \R^3$ such that $\sigma_{\Sigma^e}^{-2}\|\Psi-\textup{Id}\|_{H^2(\Sigma^e)}\leq C \|\overset{\circ}{A^e}\|_{L^2(\Sigma^e)}$ for a universal constant $C$. Using \eqref{L2norms}, \eqref{sobH2},  and Lemma \ref{metzger}, we find
$$
\|\Psi-\textup{Id}\|_\infty \leq C \sigma \|\overset{\circ}{A^e}\|_{L^2(\Sigma)} \leq c_0 \sigma^{\frac12-\delta},
$$
for some $c_0=c_0(\overline c,\delta,B_1)$. It follows
\begin{equation*}
|\vec z_\Sigma-\vec z_0|\leq \fint_\Sigma |\Psi-\textup{Id}| \ d\mu^e + c_0\sigma^{\frac12-\delta} \leq  c_0\sigma^{\frac12-\delta},
\end{equation*}
for a possibly different $c_0$ depending on the same values. \\
In addition, using the $L^\infty$ bound on $|\overset{\circ}{A}|$ from part (iii), we can apply \cite[Cor. E.1]{nerz1} which states that $\Sigma$ can be also written as a graph over the same sphere $\mathbb{S}_{\sigma_{\Sigma^e}}(\vec z_0)$, where the graph function $f$ satisfies  $\|f\|_{W^{2,\infty}(\mathbb{S}_{\sigma_{\Sigma^e}})}\leq c\sigma^{\frac12-\delta}$ for a constant $c$ depending on the same parameters as $B_\infty$. Using Lemma \ref {metzger} we see that $|\sigma_\Sigma-\sigma_{\Sigma^e}| \leq C \sigma^{\frac 12-\delta}$ and that
$f$ satisfies an analogous $W^{2,\infty}$ estimate when considered as a map on $\mathbb{S}_{\sigma_{\Sigma}}(\vec z_0)$ with the metric $g$. \medskip\\
(v) The estimate follows from the $p>2$ generalization of De Lellis-M\"uller's estimate due to Perez \cite[Thm. 1.1]{perez}, with a remainder term coming from the Riemannian asymptotically flat metric.
\end{proof}

\section{Spectral theory}

\subsection{Mass and the stability operator}
In this section we study the spectral properties of the stability operator associated to a round surface and collect other auxiliary results that will be needed in the analysis of the flow afterwards.

Unless explicitly stated, in this section we  consider closed surfaces $\Sigma$ belonging to a roundness class  $\mathcal{W}^\eta_\sigma(B_1,B_2)$
for fixed parameters $\eta,B_1,B_2$ and a general large $\sigma$. It is tacitly meant that the constants $c$ and $\sigma_0$ which appear in the statements only depend on $\eta,B_1,B_2$ and on the constants $\bar c,\delta$ in Definition \ref{asymflat}.

We begin by recalling the definition of Hawking mass of a surface.

\begin{dfn} Let $(\M,\overline \g)$ be a 3-dimensional manifold, and $\iota:\Sigma\hookrightarrow \M$ be a closed surface.  The \emph{Hawking mass} of $\Sigma$ is defined as
\begin{equation}\label{hawmass311}
m_H(\Sigma):=\sqrt{\frac{|\Sigma|}{16\pi}}\left(1-\frac{1}{16\pi}\int_\Sigma H^2 \ d\mu\right).
\end{equation}
\end{dfn}

Using the equivalent definition of ADM-mass in terms of the Einstein tensor \cite{miaotam} and the Gauss-Bonnet theorem one can show that the Hawking mass of a round surface $\Sigma$ is asymptotic to the ADM-mass of $\M$ for large radius. More precisely, the results of \cite[Appendix A]{nerz1} give the following:

\begin{prp}\label{equiv-mass}
There exist $c$ and $\sigma_0$ such that, for any $\Sigma \in \mathcal{W}^\eta_\sigma(B_1,B_2)$ with $\sigma \geq \sigma_0$, we have
$$
|m_H(\Sigma)-\overline{m}_\textup{ADM}| \leq c \sigma^{-\delta}.
$$
\end{prp}

We now introduce the stability operator, which occurs as the second variation of the area functional.

\begin{dfn} Given a surface $\iota:\Sigma\hookrightarrow \M$ and $\textup f\in H^2(\Sigma)$,  the \emph{stability operator} associated to $\Sigma$, $L^\Sigma:H^2(\Sigma)\to L^2(\Sigma)$, is defined as 
\begin{equation*}
L^\Sigma  f:=-\Delta^\Sigma f-(|A|^2+\overline{\textup{Ric}}(\nu,\nu)) f,
\end{equation*}
where $\Delta^\Sigma$ is the Laplace-Beltrami operator on $\Sigma$.
We simply write $\Delta, L$ instead of $\Delta^\Sigma, L^\Sigma$ whenever the choice of the surface $\Sigma$ is not ambiguous. 
\end{dfn}

In \cite{nerz1} and \cite{cederbaum} the spectral properties of $L^\Sigma$ are studied under the assumption that $\Sigma$ is a round surface with constant mean curvature, showing that $L^\Sigma$ is invertible if $\overline{m}_\textup{ADM} \neq 0$ and positive definite on functions with zero mean if $\overline{m}_\textup{ADM} > 0$. Here we generalize this analysis to round surfaces where we only assume that $H$ has a small oscillation as in \eqref{cond3defround}. We will see that the positivity property of $L^\Sigma$ when $\overline{m}_\textup{ADM} > 0$ is no longer true, but that the error terms admit an estimate that will be enough for our purposes.

We first recall some properties of the Laplace-Beltrami operator on a round sphere $\mathbb{S}_\sigma (0)\subset \R^3$ with the Euclidean metric. 
On a general closed surface, the eigenvalues of $\Delta$ are all positive, except the first one which is zero, with eigenspace given by the constant functions. For the Euclidean sphere, the first nonzero eigenvalue has multiplicity three and is given by
$$
\lambda^e_\alpha=\frac{2}{\sigma^2}, \ \alpha=1,2,3.
$$
An orthonormal basis for the eigenspace is given by the normalized coordinate functions restricted on $\mathbb{S}_\sigma (0)$
$$
f_\alpha^e(x)  =\sqrt{\frac{3}{4\pi \sigma^4}} x_\alpha, \ \alpha=1,2,3.
$$
The remaining eigenvalues satisfy the bound
$$
\lambda^e_i \geq \lambda^e_4=\frac{6}{\sigma^2}, \ \forall \, i \geq 4.
$$

Given a round surface $\Sigma \in \mathcal{W}^\eta_\sigma(B_1,B_2)$, we know from Lemma \ref{cor1}(iv) that $\Sigma$ can be written as a graph over a Euclidean sphere. Similarly to \cite{cederbaum}, we can use this map to identify functions on $\Sigma$ with functions on $\mathbb{S}_{\sigma_{\Sigma}}(0)$. We recall the statement of Lemma 2 of \cite{cederbaum}, which measures how much the first eigenvalues and the corresponding eigenfunctions of the Laplace-Beltrami operator on a round surface in the physical metric differ from the ones of the approximating sphere in the Euclidean metric. 

\begin{lem}
\label{ortonormalsystemL2} 
There exist $c,\sigma_0>0$ such that, for any $\Sigma \in \mathcal{W}^\eta_\sigma(B_1,B_2)$ with $\sigma \geq \sigma_0$, there is a complete orthonormal system in $L^2(\Sigma)$ consisting of eigenfunctions $\{f_\alpha\}_{\alpha=0}^\infty$ such that 
\begin{equation*}
-\Delta f_\alpha=\lambda_\alpha f_\alpha, \qquad \text{with $0=\lambda_0<\lambda_1\leq \lambda_2\leq...$}
\end{equation*}
Moreover, after possibly a rotation in the Euclidean coordinates we have,
for $\alpha=1,2,3$ ,
\begin{equation*}
\left|\lambda_\alpha-\frac{2}{\sigma_\Sigma^2}\right|\leq c\sigma^{-\frac52-\delta}, \quad \|f_\alpha-f_\alpha^e\|_{W^{2,2}(\Sigma)}\leq c\sigma^{-\frac12-\delta},
\end{equation*}
where $f_\alpha^e=\sqrt{\frac{3}{4\pi \sigma_\Sigma^4}} x_\alpha$ defined on $\R^3$. In addition, for $\alpha,\beta=1,2,3$ we have
\begin{equation}\label{ineq12stim}
\int_\Sigma \left|\langle \nabla f_\alpha,\nabla f_\beta\rangle-\frac{3\delta_{\alpha\beta}}{\sigma_\Sigma^2|\Sigma|}+\frac{f_\alpha f_\beta}{\sigma_\Sigma^2}\right| \ d\mu_g\leq c\sigma^{-\frac52-\delta}.
\end{equation}
The remaining eigenvalues satisfy
\begin{equation}\label{highmagnDelta}
\lambda_\alpha>\frac{5}{\sigma_\Sigma^2}, \quad \forall  \alpha>3.
\end{equation}
\end{lem}

We observe that in \cite{cederbaum} the above statement is given under the additional requirement that $\Sigma$ has constant space-time mean curvature; however, the proof works in the same way under the assumption that $H$ has a small oscillation as in \eqref{LinfcontrolH}.

As in \cite{cederbaum,nerz1}, we introduce a spectral decomposition for functions on $\Sigma$. We write $\langle u,v \rangle_2$ for the $L^2(\Sigma)$ scalar product of functions $u,v:\Sigma \to \R$.

\begin{dfn}\label{spectral} Given a surface $\Sigma \in \mathcal{W}^\eta_\sigma(B_1,B_2)$, let $\{ f_\alpha \}$, $\alpha=1,2,3$, be as in the previous lemma. For every $\textup w\in L^2(\Sigma)$ we define
\begin{equation}\label{spectraldec}
\textup w^t:=\sum_{\alpha=1}^3\langle \textup w,f_\alpha\rangle_2f_\alpha, \qquad \textup w^d:=\textup w-\textup w^t.
\end{equation}
We call 
$\textup w^t$ the \emph{translational part} and $\textup w^d$ the \emph{difference part} of $\textup w$.
\end{dfn}

To proceed in the analysis of the stability operator $L$, the next step is to estimate the contribution of the term containing $\overline{\textup{Ric}}(\nu,\nu)$.

\begin{prp}\label{prop551}
There exist $c,\sigma_0>0$ such that, for any $\Sigma \in \mathcal{W}^\eta_\sigma(B_1,B_2)$ with $\sigma \geq \sigma_0$ we have for every $\alpha,\beta\in\{1,2,3\},$ with $\alpha\neq \beta$,
\begin{equation*}
\left|\int_\Sigma \left(\overline{\textup{Ric}}(\nu,\nu)-\frac{H^2-h^2}{4}\right) f_\alpha f_\beta \ d\mu \right|\leq c\sigma^{-3-\delta},
\end{equation*}
and for every $\alpha\in\{1,2,3\}$ 
\begin{equation*}
\left| \lambda_\alpha- \frac{h^2}{2}-\frac{6m_H(\Sigma)}{\sigma_\Sigma^3}-\int_\Sigma \left(\overline{\textup{Ric}}(\nu,\nu)-\frac{H^2-h^2}{4}\right) f_\alpha^2  \ d\mu \right| \leq \textup c\sigma^{-3-\delta}.
\end{equation*}
\end{prp}
\begin{proof}
The proof follows the same strategy as in \cite[Lemma 4.5]{nerz1} and \cite[Lemma 3]{cederbaum}, so we will only highlight the differences due to the fact that our $\Sigma$ has not constant mean curvature. We recall that, by Lemma \ref{cor1}, we have the bounds
\begin{equation}\label{recalldecay}
H-h=O(\sigma^{-\frac{3}{2}-\delta}), \qquad 
H^2-h^2=O(\sigma^{-\frac{5}{2}-\delta}), \qquad
\|\overset{\circ}{A}\|^2 = O(\sigma^{-3-2\delta}),
\end{equation}
where the $O(\sigma^\alpha)$ notation means that the quantity is bounded in absolute value by $c \sigma^\alpha$, with $c$ depending on the usual parameters described at the beginning of the section.

By an application of the Bochner-Lichnerowicz formula we obtain, as in formula (41) in \cite{cederbaum}, the estimate
$$
\left| \lambda_\alpha^2 \delta_{\alpha \beta} - 
\int_\Sigma S \, \langle \nabla f_\alpha, \nabla f_\beta  \rangle \ d\mu \right|\leq c\sigma^{-5-\delta}.
$$
Using the Gauss equations and the bound in \eqref{recalldecay} on  $\|\overset{\circ}{A}\|^2$, we deduce
$$
\left| \lambda_\alpha^2 \delta_{\alpha \beta} - 
\int_\Sigma  \left(\overline{\textup{S}} -2 \overline{\textup{Ric}}(\nu,\nu)+\frac{H^2}{2}\right) \, \langle \nabla f_\alpha, \nabla f_\beta  \rangle \ d\mu \right|\leq c\sigma^{-5-\delta}.
$$
By writing $H^2=h^2+(H^2-h^2)$ and using the properties of the eigenfunctions, we find
$$
\left|
\left(\lambda_\alpha^2  -\lambda_\alpha \frac{h^2}2 \right)\delta_{\alpha \beta} -
\int_\Sigma \left(\overline{\textup{S}} - 2\overline{\textup{Ric}}(\nu,\nu)+\frac{H^2-h^2}{2}\right) \, \langle \nabla f_\alpha, \nabla f_\beta  \rangle \  d\mu \right|\leq c\sigma^{-5-\delta},
$$
which in view of \eqref{boundonriem}, \eqref{ineq12stim} and \eqref{recalldecay} implies
\begin{equation}\label{stima1}
\left|
\left(\lambda_\alpha^2  -\lambda_\alpha \frac{h^2}2 \right)\delta_{\alpha \beta} -
\int_\Sigma \left(\overline{\textup{S}} - 2\overline{\textup{Ric}}(\nu,\nu)+\frac{H^2-h^2}{2}\right) \, \left( \frac{3\delta_{\alpha\beta}}{\sigma_\Sigma^2|\Sigma|}-\frac{f_\alpha f_\beta}{\sigma_\Sigma^2} \right) \, d\mu \right|\leq c\sigma^{-5-\delta}.
\end{equation}
Compared to the proof in \cite{cederbaum,nerz1}, we have the additional term with $H^2-h^2$. One part of this contribution can be estimated observing that
 \begin{equation}
\int_\Sigma \frac{H^2-h^2}{2}\frac{3\delta_{\alpha\beta}}{\sigma_\Sigma^2|\Sigma|}\ d\mu =
\frac{3\delta_{\alpha\beta}}{2 \sigma_\Sigma^2|\Sigma|} \int_\Sigma (H-h)^2 d\mu
=  O(\sigma^{-5-2\delta}), \label{stima2}
  \end{equation}
  since $\int(H^2-h^2) d\mu=  \int(H-h)^2 d\mu$. The term containing $(H^2-h^2)f_\alpha f_\beta$, on the contrary, cannot be absorbed in the $O(\sigma^{-5-\delta})$ remainder and will be left as it is.
 The remaining terms in \eqref{stima1} can be rewritten in the following way using Gauss-Bonnet theorem and the definition of $m_H$, as shown in \cite{cederbaum,nerz1},
\begin{eqnarray}
& & \left(\lambda_\alpha^2  -\lambda_\alpha \frac{h^2}2 \right)\delta_{\alpha \beta} -
\int_\Sigma \left(\overline{\textup{S}} - 2\overline{\textup{Ric}}(\nu,\nu) \right) \,  \left( \frac{3\delta_{\alpha\beta}}{\sigma_\Sigma^2|\Sigma|}-\frac{f_\alpha f_\beta}{\sigma_\Sigma^2} \right) \, d\mu
\nonumber\\
& = & 
\frac{2}{\sigma_\Sigma^2} \left(\lambda_\alpha  -  \frac{h^2}2  \right)\delta_{\alpha \beta} 
  -\frac{12m_H(\Sigma)}{\sigma_\Sigma^5} \delta_{\alpha \beta} 
- \frac{2}{\sigma_\Sigma^2} \int_\Sigma \overline{\textup{Ric}}(\nu,\nu)  f_\alpha  f_\beta   \ d\mu 
+ O(\sigma^{-5-\delta}). \label{stima3}
\end{eqnarray}
Combining formulas \eqref{stima1}, \eqref{stima2}, \eqref{stima3}
and simplifying the factor $2/\sigma_\Sigma^2$ we conclude
$$
\left|
\left(\lambda_\alpha  -  \frac{h^2}2 -\frac{6m_H(\Sigma)}{\sigma_\Sigma^3} \right )\delta_{\alpha \beta} -
\int_\Sigma \left(\overline{\textup{Ric}}(\nu,\nu)-\frac{H^2-h^2}{4}\right) \,  f_\alpha f_\beta   \ d\mu \right|\leq c\sigma^{-3-\delta},
$$
which yields the thesis.
\end{proof}

We can now describe the behaviour of the bilinear form associated to the stability operator. 

\begin{prp}\label{lem_stabilityoperonproj} 
There exist $c,\sigma_0>0$ such that any surface $\Sigma \in \mathcal{W}^\eta_\sigma(B_1,B_2)$ with $\sigma \geq \sigma_0$ satisfies the following: for any $u,\phi,w \in H^2(\Sigma)$ with 
$u \in\textup{span}\{f_1,f_2,f_3\}$,  $\varphi\in \left(\textup{span}\{f_0,f_1,f_2,f_3\}\right)^\perp$ and $w$ with zero mean value, we have
\begin{equation}\label{eq3342}
\left|\langle Lu,u\rangle_2-\frac{6m_H(\Sigma)}{\sigma_\Sigma^3}\|u\|^2_2 +\frac{3h}{2}\int_\Sigma(H-h)u^2 \, d\mu \right|\leq c\sigma^{-3-2\delta}\|u\|^2_2,
\end{equation}
\begin{equation}\label{eq3352}
\langle L\varphi,\varphi\rangle_2>\frac{2}{\sigma_\Sigma^2}\|\varphi\|_2^2,
\end{equation}
 \begin{equation}\label{coroll551}
\left| \langle Lu,w \rangle_2 \right| \leq \textup c\sigma^{-\frac52-\delta}\|u\|_2\|w\|_2.
\end{equation}
\end{prp}
\begin{proof}
Let $u,v\in \textup{span}\{f_1,f_2,f_3\}$. Using the bilinearity of $L$ and Proposition \ref{prop551} we find
\begin{eqnarray*}
\langle Lu,v\rangle_2 & = &
\left(\frac{h^2}{2}+\frac{6m_H(\Sigma)}{\sigma_\Sigma^3} \right) \langle u,v\rangle_2-\int_\Sigma \left( |A|^2 + \frac{H^2-h^2}{4} \right) uv \, d\mu   \\
& & +\, \textup{O}(\sigma^{-3-2\delta})\|u\|_2\|v\|_2 \\
& = & \frac{6m_H(\Sigma)}{\sigma_\Sigma^3}\langle u, v \rangle_2-\int_\Sigma \left(|\overset{\circ}{A}|^2+ \frac{3(H^2-h^2)}{4} \right) uv \, d\mu+\textup{O}(\sigma^{-3-2\delta}) \|u\|_2\|v\|_2.
\end{eqnarray*}
Since by \eqref{recalldecay} we have $|\overset{\circ}{A}|^2=\textup O(\sigma^{-3-2\delta})$ and $H^2-h^2=2h(H-h)+\textup O(\sigma^{-3-2\delta})$, we find
\begin{equation}\label{bilinear}
\langle Lu,v\rangle_2=\frac{6m_H(\Sigma)}{\sigma_\Sigma^3}\langle u, v \rangle_2-\frac{3h}{2}\int_\Sigma(H-h)uv \, d\mu+\textup O(\sigma^{-3-2\delta})\|u\|_2\|v\|_2.
\end{equation}
By choosing $u=v$, we obtain \eqref{eq3342}.

We now recall that, by Proposition \ref{equiv-mass}, the Hawking mass $m_H(\Sigma)$ is bounded uniformly in $\Sigma$ for $\sigma$ large. Since $\|H-h\|_{L^\infty(\Sigma)}=\textup O(\sigma^{-\frac32-\delta})$, we see that \eqref{bilinear} implies 
\begin{equation}\label{bilinear2}
 \left| \langle Lu,v \rangle_2 \right| \leq  \textup c\sigma^{-\frac52-\delta}\|u\|_2\|v\|_2.
 \end{equation}

Suppose now that $\varphi\in \left(\textup{span}\{f_0,f_1,f_2,f_3\}\right)^\perp$. Then
\begin{equation*}
\langle L\varphi,\varphi\rangle_2\geq \left(\lambda_4- \sup_\Sigma\,\left|  |A|^2+\overline{\textup{Ric}}(\nu,\nu)\right|\right)\int_\Sigma \varphi^2 \ d\mu,
\end{equation*}
by the characterization of $\lambda_4$. We have $\lambda_4 > \frac{5}{\sigma_\Sigma^2}$ by Lemma \ref{ortonormalsystemL2} and $\left| \, |A|^2+\overline{\textup{Ric}}(\nu,\nu)\right|\leq \frac{3}{\sigma_\Sigma^2}$ for $\sigma$ large by \eqref{controlon|A|}, \eqref{radiicondition38} and \eqref{boundonriem} and so we obtain \eqref{eq3352}.

To prove the last assertion, we observe that, since $w$ has zero mean value, we can write $w=w_1+w_2$ with $w_1 \in\textup{span}\{f_1,f_2,f_3\}$ and $w_2 \in \left(\textup{span}\{f_0,f_1,f_2,f_3\}\right)^\perp$. Then we have  $\langle u,w_2 \rangle_2=\langle  \Delta u,w_2 \rangle_2=0$, which implies
 \begin{equation*}
\langle Lu,w_2 \rangle_2=\int_\Sigma \left(-|A|^2-\overline{\textup{Ric}}(\nu,\nu)\right) u w_2 \, d\mu=\int_\Sigma \left(-|A|^2-\overline{\textup{Ric}}(\nu,\nu)+\frac{h^2}{2}\right) u w_2 d\mu.
\end{equation*}
Since $|A|^2-\frac{h^2}{2}=|\overset{\circ}{A}|^2+\frac{H^2-h^2}{2}, \medskip$ with  $|\overset{\circ}{A}|^2=\textup O(\sigma^{-3-\delta})$ and $H^2-h^2=\textup O(\sigma^{-\frac52-\delta})$, we obtain
$$
\left| \langle Lu,w_2 \rangle_2 \right| \leq \textup c\sigma^{-\frac52-\delta}\|u\|_2\|w_2\|_2.
$$
Combining this with  \eqref{bilinear2} with $v=w_1$, we obtain \eqref{coroll551}.
\end{proof}

\begin{rem}
We point out that  the last term in the left-hand side in \eqref{eq3342} is in general of order $ O(\sigma^{-\frac 52 -\delta}) ||u||_2^2$ and  cannot be absorbed by the term with the mass which is of order $ O(\sigma^{-3}) ||u||_2^2$. Therefore, even if we are assuming the positivity of the mass, we cannot expect that the stability operator is positive definite on the translational eigenspace for the surfaces of our class $\mathcal{W}^\eta_\sigma(B_1,B_2)$.
\end{rem}

We conclude this part with an auxiliary estimate.

\begin{lem} \label{lemaux}
There exist $c,\sigma_0>0$ such that on any surface $\Sigma \in \mathcal{W}^\eta_\sigma(B_1,B_2)$ with $\sigma \geq \sigma_0$ we have for $\alpha=1,2,3$ 
\begin{equation}\label{stability2}
\left|\left\langle L(H-h),\frac{\nu_\alpha}{\sigma}\right\rangle_{L^2(\Sigma)}\right|\leq c \sigma^{-3-2\delta}
\end{equation}
where $\nu_\alpha:=\overline\g(\nu,\overline{\textup e}_\alpha)$, with $\{ \overline{\textup e}_\alpha \}_{\alpha=1,2,3}$ the canonical basis in the Euclidean coordinates.
\end{lem}
\begin{proof}
Let $f_\alpha^e$, with $\alpha=1,2,3$, be the Euclidean eigenfunctions of $-\Delta$ on the  sphere of radius $\sigma_\Sigma$ which satisfy, according to Lemma \ref{ortonormalsystemL2}, $\|f_\alpha-f_\alpha^e\|_{W^{2,2}(\Sigma)}\leq C\sigma^{-\frac12-\delta}$.
Then we have
$$
f_\alpha^e=\sqrt{\frac{3}{4\pi\sigma_\Sigma^4}}x_\alpha=\frac{\sqrt3}{\sqrt{|\Sigma|}} \nu^e_\alpha,
$$
where we have set $\nu^e_\alpha= \overline\g^e(\nu^e,\overline{\textup e}_\alpha)$. 
Using Lemma \ref{metzger} we find
$$\|\nu_\alpha-\nu^e_\alpha\|_{H^1(\Sigma)} \leq C \sigma \|\nu_\alpha-\nu^e_\alpha\|_{W^{1,\infty}(\Sigma)} 
 \leq C\sigma^{\frac12-\delta}.$$
Therefore
\begin{equation}\label{falpha}
 \left| \left| \sqrt{\frac{|\Sigma|}{3}} f_\alpha - \nu_\alpha \right| \right|_{H^1(\Sigma)}  \leq 
 \sqrt{\frac{|\Sigma|}{3}} \|f_\alpha-f_\alpha^e\|_{H^1(\Sigma)}+  \|\nu_\alpha-\nu^e_\alpha\|_{H^1(\Sigma)}
  \leq  c\sigma^{\frac12-\delta}.
\end{equation}

Since $\left| \, |A|^2+\overline{\textup{Ric}}(\nu,\nu) \right| \leq c \sigma^{-2}$, by definition of $L$ we have, for any $u,v \in H^1(\Sigma)$, 
 \medskip
\begin{equation}
\left|\left\langle Lu, v \right\rangle_{L^2}\right|  \leq   \| \nabla u \|_{L^2} \| \nabla v \|_{L^2}
+ c \sigma^{-2} \|  u \|_{L^2} \| v \|_{L^2}
\leq \sigma^{-2} \|  u \|_{H^1} \| v \|_{H^1}. \medskip
\end{equation}
Using this, together with \eqref{falpha},\eqref{coroll551} and \eqref{L2norms}, we find
\begin{eqnarray*}
\left|\left\langle L(H-h),\nu_\alpha \right\rangle_{L^2}\right| & \leq & 
 \left|\left\langle L(H-h), \sqrt{\frac{|\Sigma|}{3}} f_\alpha - \nu_\alpha \right\rangle_{L^2}\right|
+ \left|\left\langle L(H-h),\sqrt{\frac{|\Sigma|}{3}} f_\alpha \right\rangle_{L^2}\right|
\\
& \leq & c \sigma^{-2} \|  H-h \|_{H^1} \left| \left| \sqrt{\frac{|\Sigma|}{3}} f_\alpha - \nu_\alpha \right| \right|_{H^1}
+ c \sigma^{-\frac52-\delta} \, \sqrt{\frac{|\Sigma|}{3}} \, \|  H-h \|_{L^2} \\
& \leq &  c \sigma^{-2-\delta},
\end{eqnarray*}
which implies the assertion.
\end{proof}

\subsection{The translational part of the mean curvature}

We analyze now an important property of the translational part of a function (see Definition \ref{spectral}). To explain it heuristically, consider first the case of a Euclidean round sphere $\mathbb{S}_\sigma(0) \subset \mathbb{R}^3$ and take any $u \in\textup{span}\{f^e_1,f^e_2,f^e_3\}$ i.e. $u= \sum_{i=1}^3 a_\alpha x_\alpha$ for some coefficients $a_\alpha$. Then, the integral of any odd power of $u$ on $\mathbb{S}_\sigma(0)$ vanishes for symmetry reasons. This allows to obtain a strong bound on the corresponding integral when we consider a round surface in an asymptotically flat space. We focus here on the case of a third power, which is the one that we need in the sequel.

\begin{lem}\label{approxfHtlemma}
There exist $c,\sigma_0>0$ such that on any surface $\Sigma \in \mathcal{W}^\eta_\sigma(B_1,B_2)$ with $\sigma \geq \sigma_0$
we have, for any $u \in \textup{span}\{f_1,f_2,f_3\}$,
\begin{equation*}
\left| \int_{\Sigma} u^3 \ d\mu\right|\leq c\sigma^{-\frac 32-\delta}  ||u ||^3_{L^2(\Sigma)}.
\end{equation*}
\end{lem}

\begin{proof}
As before, we write $\Sigma$ as a spherical graph over a Euclidean sphere and we identify correspondingly functions on $\Sigma$ and on $\mathbb{S}_{\sigma_\Sigma}(0)$. 
By Lemma \ref{ortonormalsystemL2} and the Sobolev immersion \eqref{sobH2}, we have
\begin{equation}\label{compare}
 || f_\alpha - f_\alpha^e ||_{L^\infty(\Sigma)} \leq C \sigma^{-1}  || f_\alpha - f_\alpha^e ||_{H^2(\Sigma)} \leq c \sigma^{-\frac32-\delta}.
\end{equation}
 If we denote by $d\mu$ and $d\mu^e$ respectively the Riemannian and Euclidean measure on $\Sigma$, and by $d\mu^e_\mathbb{S}$ the Euclidean measure on $\mathbb{S}_{\sigma_\Sigma}(0)$ we have that $d\mu-d\mu^e=\textup O(\sigma^{-\frac12-\delta})d\mu$, by Lemma \ref{metzger} and also $d\mu^e-d\mu^e_\mathbb{S}=\textup O(\sigma^{-\frac12-\delta})d\mu$, by the $W^{2,\infty}$ bound on the spherical graph function in Lemma \ref{cor1}(iv). It follows that $d\mu-d\mu_e^\mathbb{S} =\textup O(\sigma^{-\frac12-\delta})d\mu$ as well.
We now define the auxiliary function
\begin{equation*}
\tilde u:=\sum_{\alpha=1}^3\langle u,f_\alpha\rangle_{L^2(\Sigma)}f_\alpha^e,
\end{equation*}
which is a combination of the Euclidean eigenfunctions $f_\alpha^e$ with the coefficients $\langle u,f_\alpha\rangle_{L^2(\Sigma)}$ occurring in the Riemannian decomposition of $u$.
Since $|\langle u,f_\alpha \rangle_{L^2(\Sigma)}| \leq \| u \|_{L^2(\Sigma)}$ and $\|f_\alpha^e \|_{L^\infty(\Sigma)} \leq c \sigma^{-1}$,
we have
$$
\|\tilde u \|_{L^\infty(\Sigma)} \leq c \sigma^{-1}\| u\|_{L^2(\Sigma)}, \qquad
\|\tilde u \|_{L^2(\Sigma)} \leq c \| u\|_{L^2(\Sigma)},
$$
$$
| u-\tilde u |  =   \left|\sum_{\alpha=1}^3\langle u,f_\alpha\rangle_{L^2(\Sigma)} (f_\alpha - f_\alpha^e) \, \right|
\leq \|u\|_{L^2(\Sigma)} \sum_{\alpha=1}^3  | f_\alpha - f_\alpha^e |.
$$
By \eqref{compare}, we deduce
$$
\|u-\tilde u\|_{L^\infty(\Sigma)}   \leq c \sigma^{-\frac32-\delta}\| u\|_{L^2(\Sigma)}.
$$
Then we can compute
\begin{eqnarray}
\left|\int_{\Sigma} \left (u^3- \tilde u^3 \right) \ d\mu \right| & \leq & \frac 32 \|u-\tilde u\|_{L^\infty(\Sigma)} \left( \| u \|^2_{L^2(\Sigma)} + \| \tilde u\|^2_{L^2(\Sigma)} \right) \nonumber
\\ & \leq & 
c \sigma^{-\frac32-\delta} \| u \|^3_{L^2(\Sigma)}.\label{ineq}
\end{eqnarray}

When considered as a function on $\mathbb{S}_{\sigma_\Sigma}(0)$, the function $\tilde u$ satisfies for symmetry reasons
$$
\int_{\mathbb{S}_{\sigma_\Sigma}} \tilde u ^3 \ d\mu^{\mathbb{S}}_e=0.
$$
It follows
\begin{eqnarray*}
\left|\int_{\Sigma}  \tilde u^3  \ d\mu \right| & = & \left|\int_{\Sigma} \tilde u^3  d\mu  -
\int_{\mathbb{S}_{\sigma_\Sigma}} \tilde u ^3 \ d\mu^{\mathbb{S}}_e \right| 
 \leq  c\sigma^{-\frac12-\delta} \int_{\Sigma} |\tilde u|^3 d\mu \\
& \leq & c \sigma^{\frac32-\delta}\|\tilde u\|^3_{L^\infty(\Sigma)} \leq c \sigma^{-\frac32-\delta} \| u\|^3_{L^2(\Sigma)},
\end{eqnarray*}
which implies the assertion, thanks to \eqref{ineq}.
\end{proof}

In particular the above estimate can be applied to the translational part of the mean curvature. Observe that by definition $H^t=(H-h)^t$; we use the longer expression $(H-h)^t$ because it makes more explicit the relation with the speed of the flow studied in the next section and the property of zero mean value.

\begin{cor}\label{approxfHtlemma2}
There exist $c,\sigma_0>0$ such that on any surface $\Sigma \in \mathcal{W}^\eta_\sigma(B_1,B_2)$ with $\sigma \geq \sigma_0$
the translational part of the mean curvature satisfies
\begin{equation*}
\left| \int_{\Sigma} ((H-h)^t)^3 \ d\mu\right|\leq c\sigma^{-2-2\delta}\|(H-h)^t\|_{L^2(\Sigma)}^2
\end{equation*}
\end{cor}
\begin{proof}
By \eqref{L2norms}, we have   $\|(H-h)^t\|_{L^2} \leq \|H-h\|_{L^2} \leq c\sigma^{-\frac 12-\delta}$. Then the assertion follows from the previous lemma.
\end{proof}
\begin{rem} Observe that this result is sharper than the one we would obtain by simply estimating one factor of $(H-h)^t$ with $\|H-h\|_{L^\infty(\Sigma)}=O(\sigma^{-\frac 32-\delta})$.
\end{rem}

The next lemma provides an approximation of $\|(H-h)^t\|_{L^2}$ which will be useful in the sequel.

\begin{lem}\label{lempit}
There exist $c,\sigma_0>0$ such that on any surface $\Sigma \in \mathcal{W}^\eta_\sigma(B_1,B_2)$ with $\sigma \geq \sigma_0$ we have,  for any $\varepsilon>0$,
\begin{eqnarray*}
\left| \frac{4\pi}{3} \|(H-h)^t\|^2_{L^2(\Sigma)} -
 \sum_{\alpha=1}^3\left\langle H-h,\frac{\nu_\alpha}{\sigma_\Sigma}\right\rangle_{L^2(\Sigma)}^2 \right|
 & \leq &  c \sigma^{-1-2\delta} (1+\varepsilon^{-1}) \|H-h \|^2_{L^2(\Sigma)}
 \\ & &
 +  \varepsilon \|(H-h)^t \|^2_{L^2(\Sigma)}.
\end{eqnarray*}
\end{lem}
\begin{proof}
We have
\begin{eqnarray*}
&& 
 \frac{4\pi}{3} \|(H-h)^t\|^2_{L^2(\Sigma)} - \sum_{\alpha=1}^3\left\langle H-h,\frac{\nu_\alpha}{\sigma_\Sigma}\right\rangle_{L^2(\Sigma)}^2\\
& = & 
\sum_{\alpha=1}^3 \left(
\frac{4\pi}{3} \left\langle H-h,f_\alpha \right\rangle_{L^2(\Sigma)}^2 
- \left\langle H-h,\frac{\nu_\alpha}{\sigma_\Sigma}\right\rangle_{L^2(\Sigma)}^2
\right) \\
& = & 
\sum_{\alpha=1}^3 
\left\langle H-h,\sqrt \frac{4\pi}{3} f_\alpha - \frac{\nu_\alpha}{\sigma_\Sigma} \right\rangle_{L^2(\Sigma)}
\left\langle H-h, \sqrt \frac{4\pi}{3} f_\alpha + \frac{\nu_\alpha}{\sigma_\Sigma} \right\rangle_{L^2(\Sigma)}.
\end{eqnarray*}
By \eqref{falpha}, we have
\begin{eqnarray*}
\left| \left\langle H-h,\sqrt \frac{4\pi}{3} f_\alpha - \frac{\nu_\alpha}{\sigma_\Sigma} \right\rangle_{L^2(\Sigma)} \right|
& \leq &  \| H- h \|_{L^2(\Sigma)} \left|\left| \sqrt \frac{4\pi}{3} f_\alpha - \frac{\nu_\alpha}{\sigma_\Sigma}\right|\right|_{L^2(\Sigma)} \\
& \leq &c \sigma^{-\frac 12-\delta} \|(H-h) \|_{L^2(\Sigma)}.
\end{eqnarray*}
We further observe
\begin{eqnarray*}
&& \left\langle H-h, \sqrt \frac{4\pi}{3} f_\alpha + \frac{\nu_\alpha}{\sigma_\Sigma} \right\rangle_{L^2(\Sigma)} \\
& = &
-\left\langle H-h,  \sqrt \frac{4\pi}{3} f_\alpha -  \frac{\nu_\alpha}{\sigma_\Sigma} \right\rangle_{L^2(\Sigma)}
+ 2 \sqrt \frac{4\pi}{3} \left\langle H-h ,  f_\alpha  \right\rangle_{L^2(\Sigma)},
\end{eqnarray*}
which implies, by definition of $(H-h)^t$,
$$
 \left| \left\langle H-h, \sqrt \frac{4\pi}{3} f_\alpha + \frac{\nu_\alpha}{\sigma_\Sigma} \right\rangle_{L^2(\Sigma)} \right| \\
\leq
c \sigma^{-\frac 12-\delta} \|(H-h) \|_{L^2(\Sigma)}
+ 2 \sqrt \frac{4\pi}{3}  || (H-h)^t||_{L^2(\Sigma)}.
$$
Putting together the above inequalities we obtain the assertion.
\end{proof}

We conclude the section by observing that the translational part of the mean curvature of a coordinate sphere satisfies an improved estimate if our ambient manifold satisfies the weak Regge-Teitelboim conditions.

\begin{lem}\label{rtspheres} Let $(\M,\overline \g)$ satisfy the $C_{1+\delta}^1$-Regge-Teitelboim conditions \eqref{wrtcond363}. Consider the Euclidean coordinate sphere $\Sigma=\mathbb{S}_r( 0)$  for some $r>1$. There exist $r_0,C>0$ depending only on $\bar c$ in \eqref{asymptc8}-\eqref{wrtcond363}, such that, if $r \geq r_0$, then $|(H-h)^t|^2_{L^2(\Sigma)} \leq C\sigma_\Sigma^{-2-2\delta}$.
\end{lem}

\begin{proof}
In this proof, we denote by $\sigma_\Sigma$ the area radius of the Euclidean sphere $\Sigma=\mathbb{S}_r( 0) $, which is asymptotic to  $r$  by Lemma \ref{metzger}. 
A direct computation shows that conditions \eqref{wrtcond363} imply, for any $\alpha=1,2,3$,
\begin{equation*}\label{nuHcont1}
|\nu_\alpha(\bar x)-\nu_\alpha(- \bar x)| \leq C r^{-1-\delta}, \quad |H (\bar x)-H(-\bar x)| \leq C r^{-2-\delta},
\end{equation*}
where $\nu_\alpha:=\overline\g(\nu,\overline{\textup e}_\alpha)$.
From this it also follows $|H (\bar x)\nu_\alpha(\bar x)+H(-\bar x)\nu_\alpha(- \bar x)| \leq C r^{-2-\delta}$. Since by \eqref{wrtcond363} $d\mu$ is antipodally symmetric on $\Sigma$
up to $O(r^{1+\delta}) d\mu$, we deduce
$$
\left| \int_{S_r(0)}  \nu_\alpha \, d\mu \right| \leq C r^{1-\delta}, \qquad \left|  \int_{S_r(0)} H \nu_\alpha \, d\mu \right| \leq C r^{-\delta},
$$
which implies
$$
\left| \int_{S_r(0)}  (H-h) \nu_\alpha \, d\mu \right| \leq \left| \int_{S_r(0)} H \nu_\alpha \, d\mu \right| + h \left| \int_{S_r(0)} \nu_\alpha \, d\mu \right| \leq C r^{-\delta}.
$$
Using again the asymptotic equivalence of $r$ and $\sigma_\Sigma$, we obtain
$$
\left\langle H-h,\frac{\nu_\alpha}{\sigma_\Sigma}\right\rangle_{L^2(\Sigma)} \leq C\sigma_\Sigma^{-1-\delta}.
$$
Since by Lemma \ref{metzger} $ \|(H-h) \|_{L^2(\Sigma)} \leq C \sigma_\Sigma^{-\frac 12 -\delta}$, we obtain from Lemma \ref{lempit} with $\varepsilon=1$
\begin{eqnarray*}
 \frac{4\pi}{3} \|(H-h)^t\|^2_{L^2(\Sigma)} 
 & \leq & C\sigma_\Sigma^{-2-2\delta}
 +  c \sigma^{-2-4\delta} +  \|(H-h)^t \|^2_{L^2(\Sigma)}.
\end{eqnarray*}
From the proof of Lemma \ref{lempit} one sees that in the case of a Euclidean coordinate sphere the constant $c$ in the above inequality only depends on $\bar c$ in \eqref{asymptc8}, and this implies the assertion. \end{proof}


\section{Volume preserving mean curvature flow}

\subsection{Definition of the flow and evolution equations}

\begin{dfn}
Let $(\M,\overline \g)$ be a Riemannian manifold, and let $\iota:\Sigma\hookrightarrow \M$ be a closed hypersurface. A time dependent family of immersions 
$F_t:\Sigma\hookrightarrow \M$, with ${t\in [0,T)}$ for some $0<T \leq +\infty$, which satisfies
\begin{equation}
\label{vpmcfsy45}
\begin{cases}
\frac{\partial}{\partial t} F_t(\cdot)=-(H(\cdot,t)-h(t))\nu(\cdot,t)\\
F_0=\iota
\end{cases}
\end{equation}
is called a solution to the \emph{volume preserving mean curvature flow}, with initial value $\iota$.
\end{dfn}

It is well-known that this flow is parabolic and it has a smooth solution at least locally in time. In the following, we always assume that the ambient manifold $(\M,\overline \g)$ is $3$-dimensional and $C_{\frac12+\delta}^2$-asymptotically flat. We write $\Sigma_t:=F_t(\Sigma)$ to denote the immersed surface at time $t$, and we call for simplicity $\Sigma_t$ the ``solution of the flow'' \eqref{vpmcfsy45} without mentioning explicitly the immersions $F_t$. We call $g(t)$ the induced metric on $\Sigma$ at time $t$ and by $d\mu_t$ the corresponding measure.

We recall the evolution equations satisfied by the main geometric quantities on $\Sigma_t$. 
We choose at each fixed time a frame $\vec{e}_\alpha$ on the ambient manifold $\M$ such that $\vec{e}_1,\vec{e}_2$  are tangent vectors on $\Sigma$ and $\vec{e}_3=\nu$.
Then the main geometric quantities on $\Sigma_t$ satisfy the following equations along the flow, see e.g. \cite{huiskenpolden}.

\begin{lem} 
\label{evolution41} Let $\{F_t\}_ {t\in [0,T)}$ be a solution of the flow \eqref{vpmcfsy45}. Then we have the equations
\begin{enumerate}[label=\textup{(\roman*)}]
\item $\frac{\partial g_{ij}}{\partial t}=-2 (H-h) h_{ij}$;
\item $\frac{\partial \nu}{\partial t}=\nabla H$;
\item $\frac{\partial}{\partial t}(d\mu_t)=-(H-h)Hd\mu_t$;
\item $\frac{\partial}{\partial t}h_{ij}=\nabla_i\nabla_j H+(H-h) \left(-h_{ik}h_j^k+\overline{\textup{Rm}}_{i3j3}\right)$;
\item $\frac{\partial H}{\partial t}=\Delta H+\textup (H-h)(|A|^2+\overline{\textup{Ric}}(\nu,\nu))$.
\end{enumerate} 
\end{lem}

Observe that the right-hand side of (v) can also be written as $-L(H-h)$, where $L$ is the stability operator associated to $\Sigma_t$. 
As an immediate consequence of the above equations we also have
\begin{equation}\label{derarea}
\frac{d}{d t} |\Sigma_t| = -\|H-h\|^2_{L^2(\Sigma_t)},
\end{equation}
\begin{equation}\label{dersecarea}
\frac{d}{d t} \|H-h\|^2_{L^2(\Sigma_t)} = -2 \langle L(H-h),H-h \rangle - \int_\Sigma H(H-h)^3 d\mu_t.
\end{equation}

We can rewrite the term $\nabla_i\nabla_j H$ in the right-hand side of (iv) by means of the Simons identity, as in Metzger \cite{metzger},
\begin{eqnarray}
\Delta h_{ij} & = & \nabla_i \nabla_j H + H h_i^l h_{lj} - |A|^2 h_{ij} + h_i^l \overline{\textup{Rm}}_{kjkl}+h^{lk}\overline{\textup{Rm}}_{lijk} \nonumber \\
& & + \nabla_j\left(\overline{\textup{Ric}}_{i\varepsilon}\nu^\varepsilon\right)+\nabla_l \left(\overline{\textup{Rm}}_{\varepsilon ijl}\nu^\varepsilon \right).
\label{simons}
\end{eqnarray}
In this form of the equality the derivatives of the Ricci and the Riemann tensor are taken with respect to the connection of $\Sigma$, in contrast to the formula used in \cite{huiskenyau}, where they are taken with respect to the ambient space. This allows to deal with these terms inside integral quantities by partial integration on $\Sigma$. Using Lemma \ref{evolution41} and Simons identity we obtain, by straightforward computations, the following result.

\begin{lem}\label{evolequations52} Along a solution of the volume preserving mean curvature flow we have
\begin{eqnarray}
\displaystyle\frac{\partial }{\partial t}|\overset{\circ}{A}|^2 & = & \Delta |\overset{\circ}{A}|^2-2|\nabla \overset{\circ}{A}|^2+\frac{2h}{H}\{|A|^4-H\textup{tr}(A^3)\}+2|A|^2\left(\frac{H-h}{H}\right)|\overset{\circ}{A}|^2 \nonumber \\
& & \label{eq46Lemma421}
+2(H-h)\overset{\circ}{h^{ij}}\overline{\textup{Rm}}_{kilj}\nu^k\nu^l-2\left(h_i^l\overline{\textup{Rm}}^k_{\,\,jkl}+h^{lk}\overline{\textup{Rm}}_{lijk}\right) h^{ij} \\
& & \nonumber 
-2\left(\nabla_j\left(\overline{\textup{Ric}}_{i\varepsilon}\nu^\varepsilon\right)+\nabla_l \left(\overline{\textup{Rm}}_{\varepsilon ijl}\nu^\varepsilon\right)\right)\overset{\circ}{{h}^{ij}}. \\
\displaystyle\frac{\partial}{\partial t}|\nabla H|^2 & = & \Delta|\nabla H|^2-2|\nabla^2H|^2+2(H-h)h^{ij}\nabla_iH\nabla_jH \nonumber \\
\label{eq99n}
& & +2(|A|^2+\overline{\textup{Ric}}(\nu,\nu))|\nabla H|^2-2\textup{Ric}^\Sigma(\nabla H,\nabla H) \\
\nonumber & & 
+2(H-h)\langle\nabla|A|^2,\nabla H\rangle+2(H-h)\langle\nabla\left(\overline{\textup{Ric}}(\nu,\nu)\right),\nabla H\rangle,
\end{eqnarray}
where $\textup{Ric}^\Sigma$ is the Ricci tensor on $\Sigma$ and $\langle\cdot,\cdot\rangle=\langle\cdot,\cdot\rangle_g$.
\end{lem}

\subsection{Evolution of integral quantities}
In this subsection we study the evolution of the integral quantities which appear in the definition of round surfaces, with the aim of studying the invariance of the class along the flow.
In our statements we will assume that our evolving surfaces satisfy properties \eqref{controlon|A|} and \eqref{radiicondition38}; in some of the results, we further require 
\begin{equation}\label{cinfBinf}
\|H-h\|_{L^\infty(\Sigma_t)}\leq c_\infty\sigma^{-\frac32-\delta}, \qquad \left\|\overset{\circ}{A}(t)\right\|_{L^\infty(\Sigma)}\leq
B_\infty\sigma^{-\frac32-\delta},
\end{equation}
for suitable $c_\infty,B_\infty$. On the other hand, we do not require a priori properties \eqref{cond2defround} and \eqref{cond3defround}. It will be important to keep explicit track of the dependence of the constants which appear in the estimates: in this way we can later fix the parameters of our roundness class in order to have invariance under the flow.

We start by estimating estimate the $L^4$ norm of $ |\overset{\circ}{A}|$. For this result, we can replace hypothesis \eqref{cinfBinf} by a milder assumption on $|H-h|$.

\begin{prp}\label{eqevApall} Let $\Sigma_t$ be a solution of the flow \eqref{vpmcfsy45} for $t \in [0,T]$ satisfying properties \eqref{controlon|A|} and \eqref{radiicondition38} for some $\sigma>1$. Suppose in addition 
\begin{equation}\label{controlon|H-h|}
\|H-h\|_{L^\infty(\Sigma_t)}\leq \displaystyle \frac{1}{20\sigma} 
\end{equation}
for all $t \in [0,T]$.
Then there exist a constant $C=C(\overline c,\delta)>0$ and a radius $\sigma_0=\sigma_0(\overline c,\delta)>0$ such that if $\sigma>\sigma_0$ then 
\begin{equation}\label{evL4_52}
\frac{\textup{d}}{\textup{d}t}\int_\Sigma |\overset{\circ}{A}|^4 \ d\mu_t\leq-2\int_\Sigma |\overset{\circ}{A}|^{2}|\nabla \overset{\circ}{A}|^2 \ d\mu_t
 -\frac{1}{2 \sigma^2} \int_\Sigma |\overset{\circ}{A}|^4 \ d\mu_t+C\sigma^{-6-4\delta}.
\end{equation} 
As a consequence, if $\displaystyle\int_\Sigma |\overset{\circ}{A}|^4 \ d\mu_0<B_1\sigma^{-4-4\delta}$ for some $B_1>2C$, then $\displaystyle\int_\Sigma |\overset{\circ}{A}|^4 \ d\mu_t<B_1\sigma^{-4-4\delta}$ for every $t\in [0,T]$.
\end{prp}

\begin{proof} From equation \eqref{eq46Lemma421} we deduce, using integration by parts
\begin{eqnarray}
\frac{\textup{d}}{\textup{d}t}\int_\Sigma |\overset{\circ}{A}|^4 \ d\mu_t & = & \int_\Sigma |\overset{\circ}{A}|^4 H(h-H) \ d\mu_t + 2 \int_\Sigma |\overset{\circ}{A}|^{2}\left(\frac{\partial}{\partial t}|\overset{\circ}{A}|^2\right) \ d\mu_t \nonumber \\
& = &  \int_\Sigma |\overset{\circ}{A}|^4 H(h-H) \ d\mu_t - 2\int_\Sigma | \nabla |\overset{\circ}{A}|^{2}  |^2 \ d\mu_t
-4\int_\Sigma |\overset{\circ}{A}|^{2}|\nabla \overset{\circ}{A}|^2 \ d\mu_t \label{prima}
 \\
& & + 4 \int_\Sigma  |\overset{\circ}{A}|^{2} \frac{h}{H}(|A|^4-H\textup{tr}(A^3)) \, d\mu_t + 
4 \int_\Sigma |A|^2\left(1-\frac{h}{H}\right)|\overset{\circ}{A}|^4 \, d\mu_t \label{seconda} \\
& & + 4 \int_\Sigma (H-h) |\overset{\circ}{A}|^2 \overset{\circ}{h^{ij}}\overline{\textup{Rm}}_{kilj}\nu^k\nu^l \, d\mu_t 
-8 \int_\Sigma \overline{\textup{Rm}}_{1212} |\overset{\circ}{A}|^4 \, d\mu_t  \label{terza} \\
& & - 4 \int_\Sigma \left(\nabla_j\left(\overline{\textup{Ric}}_{i\varepsilon}\nu^\varepsilon\right)+\nabla_l \left(\overline{\textup{Rm}}_{\varepsilon ijl}\nu^\varepsilon\right)\right)\overset{\circ}{{h}^{ij}} |\overset{\circ}{A}|^2\, d\mu_t, \label{quarta}
\end{eqnarray}
where we have used the identity
$$
-2\left(h_i^l\overline{\textup{Rm}}_{kjkl}+h^{lk}\overline{\textup{Rm}}_{lijk}\right)h^{ij}=-4|\overset{\circ}{A}|^2\overline{\textup{Rm}}_{1212},
$$
which follows from the symmetries of the Riemann tensor.

To estimate the above terms, we first observe that \eqref{controlon|A|} and \eqref{controlon|H-h|} imply
$$
\frac 1\sigma \leq  H \leq  \frac{\sqrt 5}{\sigma}, \qquad \left| 1 - \frac hH \right| \leq \frac{1}{20}, \qquad H(h-H) \leq \frac{1}{4 \sigma^2}.
$$
In addition, we recall the identity
\begin{equation*}
|A|^4-H\textup{tr}(A^3)=-2\kappa_1\kappa_2|\overset{\circ}{A}|^2.
\end{equation*}
Using again \eqref{controlon|A|}, we can estimate the positive terms in lines \eqref{prima}-\eqref{seconda} as follows
\begin{eqnarray}
\int_\Sigma |\overset{\circ}{A}|^4 H(h-H) \ d\mu_t  & + & 4 \int_\Sigma  |\overset{\circ}{A}|^{2} \frac{h}{H}(|A|^4-H\textup{tr}(A^3)) \, d\mu_t + 
4 \int_\Sigma |A|^2\left(1-\frac{h}{H}\right)|\overset{\circ}{A}|^4 \, d\mu_t \nonumber \\
& \leq & \left( \frac 14 - \frac{19}{10} + \frac 12 \right) \frac{1}{\sigma^2}  \int_\Sigma |\overset{\circ}{A}|^4  \, d\mu_t 
  \leq \frac{1}{\sigma^2}    \int_\Sigma |\overset{\circ}{A}|^4  \, d\mu_t .
\end{eqnarray}

We now consider the contribution of \eqref{terza}. Using \eqref{boundonriem} we find, for any $a>0$,
\begin{eqnarray*}
&&4 \int_\Sigma (H-h) |\overset{\circ}{A}|^2 \overset{\circ}{h^{ij}}\overline{\textup{Rm}}_{kilj}\nu^k\nu^l \, d\mu_t 
-8 \int_\Sigma \overline{\textup{Rm}}_{1212} |\overset{\circ}{A}|^4 \, d\mu_t \\
 & \leq & C \int_\Sigma \sigma^{-\frac72-\delta} |\overset{\circ}{A}|^{3} \ d\mu_t
+ C \int_\Sigma \sigma^{-\frac52-\delta} |\overset{\circ}{A}|^{4} \ d\mu_t \\
& \leq & C  \int_\Sigma \left( \left( \frac{a}{\sigma^2} + \sigma^{-\frac52-\delta} \right) |\overset{\circ}{A}|^4 
+\frac{1}{a^3} \sigma^{-8-4 \delta} \right) \ d\mu_t,
\end{eqnarray*}
where we have used Young's inequality
 and $C$ denotes as usual a constant which can change from line to line, but only depends on $\bar c, \delta$. 

To estimate the term in \eqref{quarta}, we use integration by parts and find
\begin{eqnarray*}
& &  -4 \int_\Sigma \left(\nabla_j\left(\overline{\textup{Ric}}_{i\varepsilon}\nu^\varepsilon\right)+\nabla_l \left(\overline{\textup{Rm}}_{\varepsilon ijl}\nu^\varepsilon\right)\right)\overset{\circ}{{h}^{ij}} |\overset{\circ}{A}|^2\, d\mu_t \\
& = & 4 \int_\Sigma \left( \overline{\textup{Ric}}_{i\varepsilon}\nu^\varepsilon
\nabla_j  ( \overset{\circ}{{h}^{ij}}|\overset{\circ}{A}|^{2} )
+ \overline{\textup{Rm}}_{\varepsilon ijl}\nu^\varepsilon
\nabla_l (\overset{\circ}{{h}^{ij}}|\overset{\circ}{A}|^{2} ) \, \right) \, d\mu_t \\
& \leq & C\int_\Sigma |\overline{\textup{Rm}}||\nabla\overset{\circ}{A}||\overset{\circ}{A}|^{2} \ d\mu_t
 \leq C\int_\Sigma\sigma^{-\frac 52 -\delta }|\nabla\overset{\circ}{A}||\overset{\circ}{A}|^{2} \ d\mu_t \\
 & \leq & 2 \int_\Sigma |\nabla \overset{\circ}{A}|^2|\overset{\circ}{A}|^{2} \ d\mu_t+ C \sigma^{-5-2\delta}\int_\Sigma |\overset{\circ}{A}|^{2} \ d\mu_t \\
 & \leq & 2 \int_\Sigma |\nabla \overset{\circ}{A}|^2|\overset{\circ}{A}|^{2} \ d\mu_t+ C \int_\Sigma \left( a \sigma^{-2} |\overset{\circ}{A}|^{4}+ \frac 1a \sigma^{-8-4\delta} \right) \, d\mu_t .
\end{eqnarray*}

The $|\nabla \overset{\circ}{A}|^2|\overset{\circ}{A}|^{2}$ term can be absorbed by the corresponding negative term in \eqref{prima}. Therefore, by choosing $a$ suitably small and $\sigma$ large (both depending only on $\bar c, \delta$) we conclude

\begin{equation}\label{eq458}
\frac{\textup{d}}{\textup{d}t}\|\overset{\circ}{A}\|_{L^4(\Sigma,\mu_t)}^4\leq-2\int_\Sigma |\overset{\circ}{A}|^{2}|\nabla \overset{\circ}{A}|^2 \ d\mu_t
 -\frac{1}{2\sigma^2}\|\overset{\circ}{A}\|_{L^4(\Sigma,\mu_t)}^4+C\sigma^{-6-4\delta}
\end{equation}
The last claim in our statement follows by a standard ODE comparison argument.
\end{proof}

We next estimate the rate of change of the volume preserving term $h(t)$ and of the $L^4$ norm of $H-h$.
\begin{lem} \label{rateh(t)}
Let $\Sigma_t$ be a solution of the flow \eqref{vpmcfsy45} for $t \in [0,T]$, satisfying properties \eqref{controlon|A|},  \eqref{radiicondition38} and \eqref{cinfBinf}.
Then there exist a constant $c=c(c_\infty,\overline c)>0$ and a universal constant $C_1>0$ such that 
\begin{equation}
\label{derivativeh(t)}
|\dot h(t)|\leq c\sigma^{-4-2\delta},
\end{equation}
\begin{equation}\label{oscillation}
\frac{\textup{d}}{\textup{d}t}\int_\Sigma (H-h)^4 \ d\mu_t\leq -12\int_\Sigma (H-h)^2|\nabla H|^2 \ d\mu_t+C_1 \sigma^{-2}\int_\Sigma (H-h)^4 \ d\mu_t+c\sigma^{-\frac{13}{2}-5\delta},
\end{equation}
provided $\sigma \geq \sigma_0$, for a suitable $\sigma_0=\sigma_0(\bar c, c_\infty,B_\infty)$.
\end{lem}

\begin{proof}
Similar to \cite[Lemma 14]{li} we compute,
 \begin{eqnarray*}
|\Sigma_t|\dot h(t) & = & \int_\Sigma \frac{\partial H}{\partial t} \ d\mu_t+\int_\Sigma H^2(h-H) \ d\mu_t+h\int_\Sigma (H-h)^2 \ d\mu_t \\
& = & \int_\Sigma (H-h)\left(|\overset{\circ}{A}|^2+\overline{\textup{Ric}}(\nu,\nu)\right) \ d\mu_t - \int_\Sigma (H-h)\left(\frac{H^2}{2}-Hh+h^2\right) \ d\mu_t  \\
& = & \int_\Sigma (H-h)\left(|\overset{\circ}{A}|^2+\overline{\textup{Ric}}(\nu,\nu)\right) \ d\mu_t-\frac12\int_\Sigma (H-h)^3 \ d\mu_t,
\end{eqnarray*}
using the property that $\displaystyle\int_\Sigma (H-h) \ d\mu_t=0$. By \eqref{boundonriem}, we have $\left|\overline{\textup{Ric}}(\nu,\nu)\right|\leq C\sigma^{-\frac52-\delta}$ where $C=C(\overline c)$. Then we can estimate 
\begin{equation*}
|\Sigma_t||\dot h(t)|\leq |\Sigma_t| \left( c_\infty\sigma^{-\frac32-\delta} (B_\infty^2\sigma^{-3-2\delta}+C\sigma^{-\frac52-\delta})+\frac{1}{2} c_\infty^3\sigma^{-\frac92-3\delta}\right).
\end{equation*}
After simplifying the $|\Sigma_t|$ factor, the lower order term on the right-hand side is $c_\infty C \sigma^{-4-\delta}$, and the other terms can be included in this one if $\sigma$ is large depending on $B_\infty$, $c_\infty, C$. This proves that \eqref{derivativeh(t)} holds for some $c$ only depending on $c_\infty ,\overline c$.

We now compute, using Lemma \ref{evolution41} and integration by parts,
\begin{eqnarray*}
\frac{\textup{d}}{\textup{d}t}\int_\Sigma (H-h)^4 \ d\mu_t & = &-12\int_\Sigma (H-h)^2|\nabla H|^2 \ d\mu_t+4\int_\Sigma (H-h)^4(|A|^2+\overline{\textup{Ric}}(\nu,\nu)) \ d\mu_t \\
& &
-4\dot h\int_\Sigma (H-h)^3 \ d\mu_t-\int_\Sigma H(H-h)^5 \ d\mu_t.
\end{eqnarray*}

By \eqref{controlon|A|},  and \eqref{boundonriem} we deduce that $\left| \, |A|^2+\overline{\textup{Ric}}(\nu,\nu) \right| \leq 3\sigma^{-2}$ if $\sigma$ is large enough. In addition,\eqref{cinfBinf}, \eqref{radiicondition38} and \eqref{derivativeh(t)} imply
\begin{equation}\label{eq85n}
\left|\dot h\int_\Sigma (H-h)^3 \ d\mu_t\right|\leq c(c_\infty,\overline c)\sigma^{-\frac{13}{2}-5\delta}.
\end{equation}
From this we obtain the conclusion, also observing that $|H(H-h)|\leq 5\sigma^{-2}$ in view of \eqref{controlon|A|}.
\end{proof}

We now analyze the evolution of the $L^4$ norm of $|\nabla H|$.
\begin{lem}\label{roundness45} 
Let $\Sigma_t$ be a solution to the volume preserving mean curvature flow \eqref{vpmcfsy45} for $t \in [0,T]$, which satisfies properties \eqref{controlon|A|} and\eqref{radiicondition38}. Then there exists a universal constant $C_2>0$ and a radius $\sigma_0=\sigma_0(\overline c,\delta)$ such that if $\sigma>\sigma_0$ then
\begin{equation*}
\frac{\textup{d}}{\textup{d}t}\int_\Sigma |\nabla H|^4 \ d\mu_t\leq -3\int_\Sigma |\nabla^2H||\nabla H|^2 \ d\mu_t+C_2\sigma^{-6}\int_\Sigma (H-h)^4 \ d\mu_t+C_2\sigma^{-2}\int_\Sigma |\nabla H|^4 \ d\mu_t.
\end{equation*}
\end{lem}
\begin{proof}
From \eqref{eq99n} we obtain, after integrating by parts,
\begin{eqnarray*}
\lefteqn{\frac{\textup d}{\textup dt}\int_\Sigma |\nabla H|^4 \ d\mu_t  =  \int_\Sigma |\nabla H|^4H(h-H) \ d\mu_t
- 4 \int_\Sigma \left| \nabla |\nabla H|^2 \right|^2 \ d\mu_t}  \\
& & -4\int_\Sigma |\nabla^2H|^2 |\nabla H|^2 \ d\mu_t
 +4\int_\Sigma (H-h)h^{ij}\nabla_iH\nabla_jH  |\nabla H|^2 \ d\mu_t  \\
 & &   -4 \int_\Sigma \textup{Ric}^\Sigma(\nabla H,\nabla H) |\nabla H|^2 \ d\mu_t 
\\
 & & 
-4  \int_\Sigma(H-h) \left( |A|^2+ \overline{\textup{Ric}}(\nu,\nu) \right)
\nabla\cdot\left(|\nabla H|^2\nabla H\right) 
 \ d\mu_t.
 \end{eqnarray*}
 By \eqref{controlon|A|}, we have that $H, |H-h|$ and $|A|$ are all bounded by $C\sigma^{-1}$. Using the asymptotic flatness \eqref{boundonriem} we also obtain that $| \textup{Ric}^\Sigma | \leq C\sigma^{-2}$ and $\left| |A|^2+ \overline{\textup{Ric}}(\nu,\nu)\right|  \leq C \sigma^{-2}$ for $\sigma$ enough large. Then we have, for $\varepsilon >0$ arbitrary,
 \begin{eqnarray*}
\frac{\textup d}{\textup dt}\int_\Sigma |\nabla H|^4 \ d\mu_t  
 & \leq & -4 \int_\Sigma |\nabla^2H|^2|\nabla H|^2 \ d\mu_t+C\sigma^{-2}\int_\Sigma |\nabla H|^4 \ d\mu_t
 \\
& & 
+ C\sigma^{-2} \int_\Sigma |H-h| |\nabla H|^2 |\nabla^2H|
 \ d\mu_t \\
  & \leq & \left(\frac \varepsilon2 C -4\right) \int_\Sigma |\nabla^2H|^2|\nabla H|^2 \ d\mu_t+C\sigma^{-2}\int_\Sigma |\nabla H|^4 \ d\mu_t
 \\
& & 
+ \frac C{2\varepsilon} \sigma^{-4} \int_\Sigma |H-h|^2 |\nabla H|^2
 \ d\mu_t.
\end{eqnarray*}
The assertion follows choosing $\epsilon=2/C$ and estimating the last term as follows
$$
 \sigma^{-4} \int_\Sigma |H-h|^2 |\nabla H|^2
 \ d\mu_t  \leq \sigma^{-6}\int_\Sigma (H-h)^4 \ d\mu_t+\sigma^{-2}\int_\Sigma |\nabla H|^4 \ d\mu_t.
$$
 \end{proof}

We can now estimate the weighted $W^{1,4}$-norm of $H-h$ which appears in condition \eqref{cond3defround}. We first prove separately a simple auxiliary inequality.
\begin{lem}\label{useful} Let $\Sigma\subset \M$ be a closed surface. Then we have, for every $\varepsilon>0$ and $\sigma>1$,  
\begin{equation*}
-\sigma^{-4}\int_\Sigma (H-h)^2|\nabla H|^2 \ d\mu \leq -\frac{\varepsilon}{2 \sigma^2}\int_\Sigma |\nabla H|^4 \ d\mu+\varepsilon^2\int_\Sigma |\nabla^2H|^2|\nabla H|^2 \ d\mu.
\end{equation*}
\end{lem} 
\begin{proof} Since $h$ is constant, we can write
\begin{eqnarray*}
\lefteqn{ \sigma^{-2}\int_\Sigma |\nabla H|^4 \ d\mu  =  \sigma^{-2}\int_\Sigma \langle \nabla (H-h), \nabla H \rangle \, |\nabla H|^2 \ d\mu} \\
& = &
-\sigma^{-2}\int_\Sigma (H-h)(\Delta H) |\nabla H|^2 \ d\mu-2\sigma^{-2}\int_\Sigma (H-h)g^{ij}\nabla_j Hg^{kl}\nabla_i\nabla_k H\nabla_l H \ d\mu
\\
& = &
\frac{\sqrt 2+2}{\sigma^2} \int_\Sigma |H-h| \, |\nabla^2H| |\nabla H|^2 \ d\mu  \leq 2 \int_\Sigma \left(\frac{(H-h)^2}{\varepsilon\sigma^4}+\varepsilon |\nabla^2H|^2\right) |\nabla H|^2 d\mu,
\end{eqnarray*}
which implies the assertion.
\end{proof}

\begin{lem}\label{roundnessmax14}
Let $\Sigma_t$ be a solution to the volume preserving mean curvature flow for $t \in [0,T]$, which satisfies properties \eqref{controlon|A|},  \eqref{radiicondition38} and \eqref{cinfBinf}. Suppose in addition that 
\eqref{cond2defround} holds for $t \in [0,T]$ for some $B_1>0$. For $\eta>0$, let us set
$$
\textup{a}_\eta(t):=\eta \sigma^{-4}\|H-h\|_{L^4(\Sigma_t)}^4+\|\nabla H\|_{L^4(\Sigma_t)}^4.
$$
Then we can find a universal constant $\eta_{\textup w}>0$, a constant $\tilde c= \tilde c(B_1,\overline c,\delta)$ and a radius $\sigma_0=\sigma_0(B_\infty,B_1,c_\infty,\delta,\overline c)>1$ such that, 
if $B_2>\tilde c$ and  $\sigma>\sigma_0$ we have the implication
$$\textup a_{\eta_w}(0)<B_2\sigma^{-8-4\delta} \  \Longrightarrow \ \textup a_{\eta_w}(t)<B_2\sigma^{-8-4\delta} \mbox{ for every }t\in[0,T].
$$
\end{lem}
\begin{proof}
From the previous Lemmas, we have that
\begin{eqnarray*}
\dot{\textup a}_\eta (t) & \leq & 
-3\int_\Sigma |\nabla^2H||\nabla H|^2 \ d\mu_t+C_2 \sigma^{-6}\int_\Sigma (H-h)^4 \ d\mu_t+C_2\sigma^{-2}\int_\Sigma |\nabla H|^4 \ d\mu_t+
\\
& & 
-12 \eta \sigma^{-4}\int_\Sigma(H-h)^2|\nabla H|^2 \ d\mu_t+C_1 \eta \sigma^{-6}\int_\Sigma (H-h)^4 \ d\mu_t+ \eta c(c_\infty,\overline c)\sigma^{-\frac{21}{2}-5\delta} \\
& 
\leq & \left(12 \eta \varepsilon^2-3\right)\int_\Sigma |\nabla^2H|^2|\nabla H|^2 \ d\mu_t +  (C_2- 6 \eta \varepsilon )\sigma^{-2}\int_\Sigma |\nabla H|^4 \,  d\mu_t \\
& &
+(\eta C_1+C_2)\sigma^{-6}\int_\Sigma (H-h)^4 \ d\mu_t+\eta c(c_\infty,\overline c)\sigma^{-\frac{21}{2}-5\delta},
\end{eqnarray*}
with $C_1,C_2$ universal constants. 
If we now choose
\begin{equation}\label{choice}
\eta=\eta_{\textup w}:=\frac 49 C^2_2, \qquad \varepsilon=\frac 34  C_2
\end{equation}
the inequality becomes
\begin{eqnarray*}
\dot{\textup a}_{\eta}(t)  & \leq &  -C_2 \sigma^{-2}\int_\Sigma |\nabla H|^4+(\eta C_1+C_2) \sigma^{-6}\int_\Sigma (H-h)^4 \ d\mu_t+\eta c(c_\infty,\overline c)\sigma^{-\frac{21}{2}-5\delta} \\
& \leq & -C_2 \sigma^{-2} a_\eta(t)+ \widetilde C \left( \sigma^{-6}\int_\Sigma (H-h)^4 \ d\mu_t+ \sigma^{-10-4\delta} \right),
\end{eqnarray*}
for another universal constant $\widetilde C>0$, where we have also used that
$$c(c_\infty,\overline c)\sigma^{-\frac{21}{2}-5\delta} \leq \sigma^{-10-4\delta}
$$
if $\sigma \geq \sigma_0$ for a suitable $\sigma_0(c_\infty,\bar c)$. 
Now we use point (v) of Lemma \ref{cor1} together with \eqref{cond2defround}, to obtain that there exists a constant $\cper=\cper(\overline c,\delta)>0$ such that 
\begin{equation*}
\int_\Sigma (H-h)^4 \ d\mu_t\leq \cper^4\left(\|\overset{\circ}{A}\|_{L^4(\Sigma,\mu_t)}^4+\sigma^{-4-4\delta}\right)\leq \cper^4(B_1^4+1)\sigma^{-4-4\delta},
\end{equation*} 
for $\sigma \geq \sigma_0$ for a suitable $\sigma_0=\sigma_0(B_\infty)$. 
We conclude that
\begin{eqnarray*}
\dot{\textup a}_\eta(t)  & \leq  & -C_2 \sigma^{-2} a_\eta(t) + \widetilde C \left( \cper^4(B_1^4+1) + 1 \right)  \sigma^{-10-4\delta} \\
& = & - C_2 \sigma^{-2} (a_\eta(t) - \tilde c  \sigma^{-8-4\delta}),
\end{eqnarray*}
for $\tilde c=\tilde c(B_1,\overline c,\delta)$. The conclusion follows by an ODE comparison argument.
\end{proof}

From now on, when considering the roundness class  $\mathcal{W}^\eta_\sigma(B_1,B_2)$, we fix the parameter $\eta$ equal to the value $\eta_w$ given by the previous Lemma, and we will no longer need to specify the dependence on $\eta$ of the constants in the estimates.

\subsection{Evolution of the barycenter and convergence}

An important assumption in the previous results was the comparability between $r_\Sigma$ and $\sigma$ in \eqref{radiicondition38}, in particular the lower bound on $r_\Sigma$ which shows that $\Sigma_t$ stays enough far from the coordinate origin to ensure the desired decay of the ambient curvature. To justify this assumption, we study now the evolution of the barycenter under the flow. We start by proving an important decay estimate on the $L^2$-norm of $H-h$, which relies on the spectral analysis of Section 3. 

\begin{prp}\label{evolutionH-h} Let $(\M,\overline \g)$ be a \textup{$C_{\frac12+\delta}^2$-asymptotically flat manifold} with  $\overline m_{\textup{ADM}}>0$. 
Given $B_1,B_2>0$, there exists $\sigma_0=\sigma_0(B_1,B_2,\overline c,\delta,\overline m_{\textup{ADM}})$ such that, if 
$\Sigma_t$ is a solution to the volume preserving mean curvature flow \eqref{vpmcfsy45} which satisfies 
$\Sigma_t\in \mathcal{W}^\eta_\sigma(B_1,B_2)$ for some $\sigma>\sigma_0$, then
\begin{equation*}
\frac{\textup d}{\textup dt}\|H-h\|_{L^2(\Sigma_t)}^2\leq -\frac{4\overline m_{\textup{ADM}}}{\sigma_{\Sigma_t}^3}\|(H-h)^t\|_{L^2(\Sigma_t)}^2-\frac{2}{\sigma_{\Sigma_t}^2}\|(H-h)^d\|_{L^2(\Sigma_t)}^2.
\end{equation*}
\end{prp}
\begin{proof} 
Let us consider the stability operator associated to $H-h$. By writing $H-h=(H-h)^t+(H-h)^d$ and applying Proposition \ref{lem_stabilityoperonproj}, we obtain
\begin{eqnarray*}
&& \langle L(H-h),H-h\rangle_2 \\
& \geq &\frac{6m_H(\Sigma_t)}{\sigma_\Sigma^3}\|(H-h)^t\|_2^2-\frac{3h}{2}\int_\Sigma (H-h)((H-h)^t)^2 \ d\mu_t -c\sigma^{-3-2\delta}\|(H-h)^t\|_2^2 \\
&& +\frac{2}{\sigma_\Sigma^2}\|(H-h)^d\|_2^2-c\sigma^{-\frac52-\delta}\|(H-h)^t\|_2\|(H-h)^d\|_2
\\
& \geq & \frac{5m_H(\Sigma_t)}{\sigma_\Sigma^3}\|(H-h)^t\|_2^2-\frac{3h}{2}\int_\Sigma (H-h)((H-h)^t)^2 \ d\mu_t+\frac{3}{2\sigma_\Sigma^2}\|(H-h)^d\|_2^2,
\end{eqnarray*}
for $\sigma$ large enough. Here we write for simplicity $\sigma_{\Sigma}$ instead of $\sigma_{\Sigma_t}$ for the area rdius of $\Sigma_t$. Since by Proposition \ref{equiv-mass} $m_H(\Sigma_t)\geq \frac{\overline m_{\textup{ADM}}}{2}>0$ for $\sigma$ sufficiently large,  we find by \eqref{dersecarea}
\begin{eqnarray*}
\frac{\textup d}{\textup{d}t}\|H-h\|_2^2 & = & -2\langle L(H-h),H-h\rangle_2-\int_\Sigma H(H-h)^3 \ d\mu_t
\\
& \leq & -\frac{5\overline m_{\textup{ADM}}}{\sigma_\Sigma^3}\|(H-h)^t\|_2^2+3h\int_\Sigma (H-h)((H-h)^t)^2 \ d\mu_t\\
& & -\frac{3}{\sigma_\Sigma^2}\|(H-h)^d\|_2^2 -h\int_\Sigma(H-h)^3 \ d\mu +c\sigma^{-3-2\delta}\int_\Sigma (H-h)^2 \ d\mu_t,
\end{eqnarray*}
where we have used $|H-h| \leq c_\infty \sigma^{-\frac 32-\delta}$ from Lemma \ref{cor1}(i). We want to show that the contribution of the $H-h$ integrals can be bounded by the remaining negative terms. By writing again $H-h=(H-h)^t+(H-h)^d$ we obtain
\begin{eqnarray*}
&& 3\int_\Sigma (H-h)((H-h)^t)^2 \ d\mu -\int_\Sigma(H-h)^3 \ d\mu \\
& = & 2 \int_\Sigma ((H-h)^t)^3 \, d\mu - \int_\Sigma 3(H-h)^t((H-h)^d)^2 \ d\mu -  \int_\Sigma ((H-h)^d)^3 \, d\mu.
\end{eqnarray*}
The first integral was considered in Corollary \ref{approxfHtlemma2}. The remaining ones can be estimated as follows
\begin{eqnarray*}
& & \left| \int_{\Sigma} 3 (H-h)^t((H-h)^d)^2 \ d\mu + ((H-h)^d)^3 \ d\mu\right|\\
& \leq & \|3(H-h)^t+(H-h)^d\|_\infty\int_\Sigma ((H-h)^d)^2 \ d\mu\leq c\sigma^{-\frac32-\delta}\int_\Sigma ((H-h)^d)^2 \ d\mu
\end{eqnarray*}
since also $\|(H-h)^d\|_{L^\infty(\Sigma)}\equiv \|(H-h)-(H-h)^t\|_{L^\infty(\Sigma)}=\textup O(\sigma^{-\frac32-\delta})$. Combining this with Corollary \ref{approxfHtlemma2}, and using $h= O(\sigma^{-1})$ we get 
\begin{eqnarray*}
\frac{\textup d}{\textup{d}t}\|H-h\|_2^2 & \leq & c \sigma^{-3-2\delta} \|(H-h)^t\|_2^2 + c \sigma^{-\frac52-\delta} \|(H-h)^d\|_2^2
\\
& & -\frac{5\overline m_{\textup{ADM}}}{\sigma_\Sigma^3}\|(H-h)^t\|_2^2-\frac{3}{\sigma_\Sigma^2}\|(H-h)^d\|_2^2
\\
 & \leq &-\frac{4\overline m_{\textup{ADM}}}{\sigma_\Sigma^3}\|(H-h)^t\|_2^2-\frac{2}{\sigma_\Sigma^2}\|(H-h)^d\|_2^2
\end{eqnarray*}
for $\sigma$ large enough.
\end{proof}

\begin{rem} An immediate consequence is that, if $\sigma$ is sufficiently large, we have 
\begin{equation}\label{estHh}
\frac{\textup{d}}{\textup dt}\|H-h\|_2^2\leq -\frac{4 \overline m_{\textup{ADM}}}{\sigma_\Sigma^3}\|H-h\|_2^2
\end{equation}
\end{rem}

Let us now define
\begin{equation*}
\Pi(t):=\sqrt{\sum_{\alpha=1}^3\left\langle H-h,\frac{\nu_\alpha}{\sigma}\right\rangle_{L^2(\Sigma_t)}^2}.
\end{equation*}
Since the $L^2$-norm of $\frac{\nu_\alpha}{\sigma}$ is uniformly bounded on a round surface for $\sigma$ large, we obtain that $\Pi(t) \leq C \|H-h\|_{L^2}$ for some universal constant $C$. We observe that, by Lemma \ref{lempit}, $\Pi(t)$ allows to estimate  the $L^2$-norm of the translational part of $H-h$ up to a constant factor. For our purposes it suffices to consider the statement of the Lemma for $\varepsilon=1$
$$
\left| \frac{4\pi}{3}\|(H-h)^t\|_{L^2}^2 - \frac{\sigma^2}{\sigma_\Sigma^2} \Pi(t)^2 \right| \leq \|(H-h)^t\|^2 + c \sigma^{-2-4\delta},
$$
where we have also used \eqref{L2norms}.
Using \eqref{radiicondition38} to obtain a rough estimate of $\sigma/\sigma_\Sigma$ we deduce that, for any $c>0$ we have, for $\sigma$ large
\begin{equation}\label{lempit2}
\|(H-h)^t\|_{L^2} \leq c \sigma^{-1-\delta} \ \Longrightarrow \Pi(t) \leq 3 c \sigma^{-1-\delta}
\end{equation}
and similarly
\begin{equation}\label{lempit3}
\Pi(t) \leq c\sigma^{-1-\delta} \ \Longrightarrow \|(H-h)^t\|_{L^2} \leq c \sigma^{-1-\delta}.
\end{equation}

We can now start our study of the barycenter of the evolving surface.

\begin{lem}\label{derivbaryc}
Let $\Sigma_t$ be a solution of the flow \eqref{vpmcfsy45} which belongs to $\mathcal{W}_\sigma^\eta(B_1,B_2)$ for all $t$. Then we have
\begin{equation}\label{speedbar}
\left| \frac{d}{dt}\vec{z}_{\Sigma_t} \right| 
\leq C \sigma^{-1}\|H-h\|_{L^2} ,
\end{equation}
\begin{equation}\label{derPi}
\left| \frac{d}{dt}\Pi(t) \right| \leq c\sigma^{-3-2\delta},
\end{equation}
for a universal constant $C$ and for a suitable $c=c(B_1,B_2)$, provided $\sigma$ is large.
\end{lem}

\begin{proof}
Let us write for simplicity $\vec z(t)=\vec z_{\Sigma_t}$.
A straightforward computation, see for example \cite{corvinowu}, shows 
\begin{equation}\label{baryev4129}
\frac{d}{dt} \vec z(t) =\frac{1}{|\Sigma_t|}\int_\Sigma (h-H)\left[\nu+ H\left(F_t(x)-\vec z(t)\right)\right] \ d\mu_t.
\end{equation}
By Lemma \ref{cor1}(iv), we have
$$
h=2\sigma_\Sigma^{-1}+ O(\sigma^{-\frac 32-\delta}), \qquad F_t(x)-\vec z(t)= \sigma_{\Sigma_t}\nu + O(\sigma^{\frac 12-\delta}). 
$$
which implies $H (F_t(x)-\vec z(t) )= 2\nu+O(\sigma^{-\frac 12-\delta})$.
Then we can estimate
\begin{eqnarray*}
\left| \int_\Sigma (h-H)\left[\nu+ H\left(F_t(x)-\vec z(t)\right)\right] \ d\mu_t \right|
 & \leq & \left| 3 \int_\Sigma (h-H)  \nu \ d\mu_t \right| \\
 & & 
+ O(\sigma^{-\frac 12-\delta})  \int_\Sigma |h-H| d\mu_t \\
& \leq & 
3  \sigma \Pi(t) +c   \sigma^{- \frac 12-\delta} |\Sigma_t|^{\frac 12}\|H-h\|_{L^2}.
\end{eqnarray*}
Since $|\Sigma_t|=O(\sigma^2)$, we see that \eqref{baryev4129} implies 
$$
\left| \frac{d}{dt}\vec{z}(t) \right| \leq C \sigma^{-1} \Pi(t) + c \sigma^{-\frac 32-\delta}\|H-h\|_{L^2},
$$
and the right-hand side is smaller than $C \sigma^{-1} \|H-h\|_{L^2}$ for $\sigma$ large, with $C$ universal constant, which yields \eqref{speedbar}. To derive \eqref{derPi} we first compute, by using Lemma \ref{evolution41},
 \begin{eqnarray*}
\frac{\textup d}{\textup dt}\left\langle H-h,\frac{\nu_\alpha}{\sigma}\right\rangle_{L^2(\Sigma)} & = & 
\frac1\sigma \int_\Sigma (-L(H-h)) \nu_\alpha  \, d\mu_t - \frac1\sigma \int_\Sigma \dot h\nu_\alpha \ d\mu_t \\
& &  
+\frac1\sigma \int_\Sigma (H-h) \langle \nabla H, e_\alpha \rangle_e \, d\mu_t- \frac1\sigma \int_\Sigma (H-h)^2H\nu_\alpha \, d\mu_t 
\end{eqnarray*}
Using \eqref{L2norms}, \eqref{stability2} and \eqref{derivativeh(t)} we conclude that the right-hand side is $O(\sigma^{-3-2\delta})$.
\end{proof}

The next result, which is similar to Proposition 3.4 in \cite{huiskenyau}, gives a bound on the possible change of area of the surface along the flow as long as it remains round.

\begin{lem}\label{HYvol}
Given $B_1,B_2$, there exist constants $c>0$ and $\sigma_0>1$ such that, if $\sigma>\sigma_0$ and $\Sigma_t$ is a solution of the flow \eqref{vpmcfsy45}  with $\Sigma_t \in  \mathcal{W}^\eta_\sigma(B_1,B_2)$ for all $t \in [0,T]$ then
$$
0 \leq \sigma_{\Sigma_0} -\sigma_{\Sigma_t} \leq c \sigma^{\frac 12-\delta}
$$ 
for every $t \in [0,T]$. 
\end{lem}
\begin{proof}
By the area-decreasing property \eqref{derarea} we have $\sigma_{\Sigma_0} \geq \sigma_{\Sigma_t}$, so we only need to prove the latter inequality. 
By \eqref{radiicondition38}, we have $r_{\Sigma_t}<\sigma/2$ and so the Euclidean coordinate sphere $\mathbb{S}_{\frac \sigma 2}(0)$ is enclosed by $\Sigma_t$ for all $t$. By definition of our flow, the compact region  between $\Sigma_t$ and $\mathbb{S}_{\frac \sigma 2}(0)$ has constant volume. We call this region $\Omega_t$. We have, using \eqref{asymflat}, 
$$
|\textup{Vol}_{e} (\Omega_t) - \textup{Vol}_{g} (\Omega_t) | \leq C \sigma^{\frac 52-\delta}.
$$
On the other hand, we know from Lemma \ref{cor1}(iv) that $\Sigma_t$ can be written as a graph over a Euclidean sphere of radius $\sigma_{\Sigma_t}$ with the radial function of order $O(\sigma^{\frac 12-\delta})$. It follows 
$$
\left| \textup{Vol}_{e} (\Omega_t) - \frac{4\pi}{3} (\sigma_{\Sigma_t}^3-r_0^3) \right|  \leq  c \sigma^{\frac 52-\delta}.
$$
Since $\textup{Vol}_{g} (\Omega_t)$ is constant, we deduce $|\sigma_{\Sigma_t}^3-\sigma_{\Sigma_0}^3 | \leq  c \sigma^{\frac 52-\delta}$, which implies the assertion.
\end{proof}

We are now ready to prove that, by an appropriate choice of the parameters of the class $\mathcal{B}^\eta_\sigma(B_1,B_2,B_\textup{cen})$ of well-centered round surfaces, see Definition \ref{roundsurf251}, and under suitable conditions on the initial surface, the solution of the flow remains inside the class for arbitrary times. Observe that, in order to control the possible drift of the barycenter, our assumptions on the initial data include an additional smallness requirement on the $L^2$-norm of the translational part of the mean curvature.

\begin{thm}\label{lemmamainth} Let $(\M,\overline \g)$ be a $C_{\frac12+\delta}^2$-asymptotically flat manifold with $\overline m_{\textup{ADM}}>0$. Let $B_1$ be chosen as in \textup{Lemma \ref{eqevApall}} and $\eta, B_2$ be chosen as in \textup{Lemma \ref{roundnessmax14}}. Then, for any $\cin>0$ there exists $\overline{B}=\overline{B}(\cin,\overline m_{\textup{ADM}})>0$ 
and $\overline{\sigma}=\overline{\sigma}(\overline c,\delta,B_1,B_2,\overline B,\overline m_{\textup{ADM}})>1$
such that the following holds. Let $F_t:\Sigma \hookrightarrow \M$ be a family of surface immersions which solve the volume preserving mean curvature flow \eqref{vpmcfsy45} for $t \in [0,T]$, and set $\Sigma_t=F_t(\Sigma)$. Suppose that the initial surface $\Sigma_0=F_0(\Sigma)$ 
satisfies, for some $\sigma \geq \overline{\sigma}$ and $B_\textup{cen} \geq \overline{B}$:
\begin{enumerate}[label=(\roman*)]
\item $\Sigma_0 \in \mathcal{B}^\eta_\sigma(B_1,B_2,B_\textup{cen})$ with $\sigma=\sigma_{\Sigma_0}$,
\item  $\|(H-h)^t\|_{L^2}\leq \cin\sigma^{-1-\delta}$ on $\Sigma_0$.
\end{enumerate}
Then $\Sigma_t \in \mathcal{B}^\eta_\sigma(B_1,B_2, 3 B_\textup{cen})$ for all $t \in [0,T]$.
\end{thm}

\begin{proof}
Let $\eta,B_1,B_2$ be fixed as in the statement and let $\cin>0$ be arbitrary. We need to prove that, if $\Sigma_0$ satisfies (i) and (ii) for suitably large $B_\textup{cen}$ and $\sigma$, whose size will be specified during the proof, then it also satisfies the conclusion.

We argue by contradiction and denote by $T'$ the first time at which the conclusion is violated, that is, $T'\in (0,T]$ is such that
\begin{equation}\label{assurdo}
(\Sigma,g(t))\in \mathcal{B}^\eta_\sigma(B_1,B_2, 3 B_\textup{cen}) \mbox{ for }t \in [0,T'), \qquad
(\Sigma,g(t))\notin \mathcal{B}^\eta_\sigma(B_1,B_2, 3 B_\textup{cen}) \mbox{ for }t =T'.
\end{equation}
This means that at least one inequality in the definition of $\mathcal{B}^\eta_\sigma(B_1,B_2,3 B_\textup{cen})$ becomes an equality at time $t=T'$. The theorem will be proved if we can exclude each of these possibilities. 

We first observe that, by Lemma \ref{HYvol}, we have $|\sigma_{\Sigma_t}-\sigma| \leq c \sigma^{\frac 12 -\delta}$ for $t \in [0,T']$. 
In addition, using the spherical graph representation in Lemma \ref{cor1} (iv), we have at any point of $\Sigma_{T'}$, 
\begin{equation*}
 \vec x  =  \vec z_0+(\sigma_\Sigma +f)\nu = \vec z(t) + \sigma \nu + O(\sigma^{\frac 12-\delta}), 
\end{equation*}
which implies, by 
\eqref{cond2defroundprima}, that 
$$
\left| \, |\vec x| - \sigma \right| \leq c \sigma^{1-\delta} 
$$
for some $c=c(B_1,B_2,B_\textup{cen},\bar c, \delta)$. From this it follows that, if $\sigma$ is large enough, the strict inequalities \eqref{radiicondition38} hold also at time $t=T'$.  Then, by parts (i)-(ii)-(iii) of Lemma \ref{cor1}, we deduce that \eqref{controlon|A|} holds as well.

Again parts (i) and (iii) of Lemma \ref{cor1} show that assumptions \eqref{cinfBinf} and \eqref{controlon|H-h|} are satisfied. Then we can apply first Lemma \ref{eqevApall} and then Lemma \ref{roundnessmax14} to prove that inequalities \eqref{cond2defround} and \eqref{cond3defround}
remain strict at time $t=T'$. We conclude that the only property of $\mathcal{B}^\eta_\sigma(B_1,B_2,3 B_\textup{cen})$ that can become an equality at time $t=T'$ is the barycenter estimate \eqref{cond2defroundprima}, and therefore we have
\begin{equation}\label{assurdo2}
|\vec z(T')| = 3 B_\textup{cen}\sigma^{1-\delta}.
\end{equation}

The rest of the proof will be devoted to exclude this last possibility, and this will require a longer work.
We first observe that, by \eqref{lempit2} and assumption (ii), we have 
\begin{equation}
\Pi(0) \leq 3 \cin \sigma^{-1-\delta}. \label{Pi-initial}
\end{equation}

Let us now denote by $t^*$ the smallest $t \in [0,T']$ such that one of the following properties holds: \\
(a) $\Pi(t) = (3 \cin +1) \sigma^{-1-\delta}$, \\
(b) $\| (H-h)^d \|_{2} \leq \| (H-h)^t \|_{2}$,  \\
(c) $|\vec{z}(t)| = 2 B_\textup{cen} \sigma^{1-\delta}$. \\
Such a $t^*$ exists, since at least property (c) must hold at some $t<T'$ because of \eqref{assurdo2}. We remark that it is possible that $t^*=0$, because property (b), in contrast with (a) and (c), is compatible with our assumptions on the initial data. 

\noindent
\textbf{Claim.} If $B_\textup{cen}$ and $\sigma$ are large enough, then at time $t=t^*$ the solution satisfies
\begin{equation}\label{baryc}
\vec{z}(t^*) < 2 B_\textup{cen} \sigma^{1-\delta},
\end{equation}
\begin{equation}\label{oscill}
\|H-h\|_{L^2(\Sigma_{t^*})} \leq c \sigma^{-1-\delta},
\end{equation}
for a suitable $c=c(B_1,B_2,\bar c, \cin)$.

If $t^*=0$ the claim is immediate. In fact, in this case the initial surface satisfies (b) and \eqref{oscill} is a consequence of hypothesis (ii) of the theorem, while \eqref{baryc} follows from (i). Therefore we assume that $t^*>0$. By definition of $t^*$, we have $\| (H-h)^d \|_{2} > \| (H-h)^t \|_{2}$ for $t \in [0,t^*)$. Then Proposition \ref{evolutionH-h} implies
\begin{equation*}
\frac{\textup d}{\textup dt}\|H-h\|_2^2\leq -\frac{2}{\sigma_{\Sigma_t}^2}\|(H-h)^d\|_2^2\leq -\frac{1}{\sigma_{\Sigma_t}^2}\|(H-h)^t\|_2^2-\frac{1}{\sigma_{\Sigma_t}^2}\|(H-h)^d\|_2^2,
\end{equation*}
for $t \in [0,t^*]$. By \eqref{radiicondition38}, we have that  $\sigma_{\Sigma_t}^2 \geq (4/5) \sigma^2$. Then we can integrate the inequality and obtain
\begin{equation}\label{estimate5307}
\|H-h\|_{L^2(\Sigma_{t})}^2\leq \|H-h\|_{L^2(\Sigma_{0})}^2e^{-\frac{4}{5\sigma^2}t}\leq  c\sigma^{-1-2\delta} e^{-\frac{4}{5\sigma^2}t}
\end{equation}
for every $t\in [0,t^*]$. Therefore, by Lemma \ref{derivbaryc},
\begin{eqnarray*}
|z(t^*)|-|z(0)| & \leq & \int_0^{t^*} c \sigma^{-1} \|H-h\|_{L^2} \, dt \leq \int_0^{t^*} c \sigma^{-\frac 32-\delta}  e^{-\frac{2}{5\sigma^2}t} \, dt \leq c \sigma^{\frac 12-\delta} .
\end{eqnarray*}
So if $B_\textup{cen}>c$ then  $|z(t^*)| < 2 B_\textup{cen} \sigma^{\frac 12-\delta}$, proving \eqref{baryc}. This also shows that case (c) of the definition of $t^*$ cannot occur, and that either (a) or (b) must hold. We  prove \eqref{oscill} dividing the two cases.
Suppose first that (a) holds. Then we have, by Lemma \ref{derivbaryc},
$$
 \sigma^{-1-\delta}  \leq  \Pi(t^*)-\Pi(0) \leq \int_0^{t^*} \left| \frac{d}{dt} \Pi(t) \right| \, dt \leq c \sigma^{-3-2\delta} t^*,
$$
which implies $t^* \geq  c^{-1} \sigma^{2+\delta}$. Substituting in \eqref{estimate5307}, we find
$$
\|H-h\|_{L^2(\Sigma_{t})}^2 \leq c\sigma^{-1-2\delta} e^{-\frac{\sigma^\delta}{c}}
$$
so that \eqref{oscill} is satisfied if $\sigma$ is large enough. If instead (b) holds at time $t=t^*$, then \eqref{oscill} follows directly by using \eqref{lempit3} to obtain 
$$
\|H-h\|_{L^2(\Sigma_{t^*})} \leq 2 \|(H-h)^t\|_{L^2(\Sigma_{t^*})}  \leq 2( 3\cin +1) \sigma^{-1-\delta}.
$$
So we have proved our claim that \eqref{baryc}-\eqref{oscill} hold at at time $t=t^*$.

To conclude the proof, we study the behaviour of the solution for $t \in [t^*,T']$. 
We can estimate $\sigma_\Sigma^3 \leq 2 \sigma^3$ using 
\eqref{radiicondition38} and deduce from \eqref{estHh}, \eqref{oscill}
\begin{equation}\label{expdecay}
\|H-h\|^2_{L^2(\Sigma_{t})} \leq \|H-h\|^2_{L^2(\Sigma_{t^*})} e^{-\frac{2 \overline m}{\sigma^3}(t-t^*)}
\leq c \sigma^{-2-2\delta} e^{-\frac{2 \overline m}{\sigma^3}(t-t^*)}, \qquad t \in [t^*,T'],
\end{equation}
where we have written for simplicity $\overline m=\overline m_{\textup{ADM}}$. 
Using Lemma \ref{derivbaryc}, we then get
\begin{eqnarray*}
|\vec z(T')-\vec z(t^*)| & \leq &  \int_{t^*}^{T'}\left| \frac{d}{dt} \vec z (t) \right| \ dt\leq 
 \int_{t^*}^{T'} c \sigma^{-1} \| H-h\|_{2} \ dt \\
 &\leq & \int_{t^*}^{T'} c \sigma^{-2-\delta} e^{-\frac{\overline m}{\sigma^3}(t-t^*)} \ dt \\
 & \leq & 
c \sigma^{-2-\delta}\left(\frac{\sigma^3}{\overline m}\right)\left(e^{-\frac{\overline m}{\sigma^3}t^*}-e^{-\frac{\overline m}{\sigma^3}T'}\right) \leq c \sigma^{1-\delta}.
\end{eqnarray*}
If $B_\textup{cen}>c$ then 
$|\vec z(T')| \leq |\vec z(t^*)| + c \sigma^{1-\delta} < (2 B_\textup{cen} +c) \sigma^{1-\delta} <3 B_\textup{cen} \sigma^{1-\delta}$, in contradiction with \eqref{assurdo2}. This shows that \eqref{assurdo2} cannot happen if $B_\textup{cen}$ and $\sigma$ are large enough, and so it concludes our proof.
\end{proof}

\begin{cor}
\label{mainthPart1} 
Let $(\M,\overline \g)$ be a $C_{\frac12+\delta}^2$-asymptotically flat manifold with $\overline m_{\textup{ADM}}>0$. For any given $c_\textup{in}>0$, let the parameters $B_1, B_2, \eta, B_\textup{cen},\sigma$ be chosen as in Theorem \ref{lemmamainth}. Then, for any initial data satisfying hypotheses (i) and (ii) of that theorem, the solution to the volume preserving mean curvature flow \eqref{vpmcfsy45} exists for every $t\in[0,\infty)$, satisfies $\Sigma_t \in \mathcal{B}^\eta_\sigma(B_1,B_2,3B_\textup{cen})$ for every $t\in[0,\infty)$ and converges exponentially in $C^\infty$ to a stable CMC surface $\Sigma_\infty \in \mathcal{B}^\eta_\sigma(B_1,B_2,3B_\textup{cen})$.  
\end{cor}

\begin{proof} Let us consider an initial surface $\Sigma_0$ which satisfies assumptions (i) and (ii) of Theorem \ref{lemmamainth}. By the local existence theory, there exists a solution of the flow defined on some maximal time interval $[0,T_{\textup{max}})$.
Then, Theorem \ref{lemmamainth} shows that $\Sigma_t$ belongs to $\mathcal{B}^\eta_\sigma(B_1,B_2,3B_\textup{cen})$ for all $t \in [0,T_{\textup{max}})$. This implies that $\Sigma_t$  is confined in a compact subset of $\M$, and therefore $\overline {\textup{Rm}}$ and each of its derivatives are uniformly bounded on $\Sigma_t$. Since the intrinsic curvature of $\Sigma_t$ is uniformly bounded by \eqref{controlon|A|}, well-known estimates based on parabolic regularity, see e.g. 
\cite{huisken86}, show that all derivatives of $A$ are also uniformly bounded. Then a standard continuation argument implies that $T_{\textup{max}}=+\infty$.

From estimate \eqref{expdecay} we see that $\|H-h\|_{L^2(\Sigma_t)}$ decays exponentially as $t \to +\infty$. Since the derivatives of any order of $H$ are uniformly bounded, interpolation estimates imply that they also decay exponentially. Then Sobolev immersion implies that $\|H-h\|_{L^\infty(\Sigma_t)}$ decays exponentially as well. Since $H-h$ is the speed of our flow, this shows that the immersions $F(\cdot,t)$ converge smoothly to a limiting map $F_\infty(\cdot)$. By standard arguments, see Lemma 8.2 in \cite{huisken1}, one can show that $F_\infty(\cdot)$ is also a smooth surface immersion and that there is exponential convergence in $C^\infty$ of the curvature. 
In particular, the limit surface $\Sigma_\infty:=F_\infty(\Sigma)$ satisfies $H \equiv h$. Then Proposition \eqref{lem_stabilityoperonproj}, 
or also Proposition 4.7 in \cite{nerz1}, show that the CMC surface $\Sigma_\infty$ is stable. Finally, the estimates in the proof of Theorem \ref{lemmamainth} for arbitrary $T>0$ show that the requirements in the definition of $\mathcal{B}^\eta_\sigma(B_1,B_2,3B_\textup{cen})$ still hold as strict inequalities on $\Sigma_\infty$.
\end{proof}

We conclude by considering the explicit example of a Euclidean coordinate sphere $\mathbb{S}_r(0)$ as initial surface for our flow. To ensure that condition (ii) of Theorem \ref{lemmamainth} is satisfied, we have to strengthen the assumptions on our ambient manifold by requiring the  $C_{1+\delta}^1$-Regge-Teitelboim conditions in Definition \ref{weakrtcond}. Then Theorem \ref{mainspheres} stated in the introduction is an immediate consequence of the previous results.

\begin{proof}[Proof of Theorem \ref{mainspheres}] By Lemma \ref{metzger}, a Euclidean sphere $\mathbb{S}_r(0)$ satisfies condition (i) of Theorem \ref{lemmamainth} if $r$ is enough large and, if $M$ satisfies the $C_{1+\delta}^1$-Regge-Teitelboim conditions, it also satisfies condition (ii) by Lemma \ref{rtspheres}. The conclusion then follows from Corollary \ref{mainthPart1}.
\end{proof}

\begin{rem}
Once the convergence of the flow is established, arguments similar to \cite{huiskenyau, nerz1} allow to conclude that the family of all limiting surfaces for large values of $r$ form a CMC-foliation of the outer part of $\M$.
\end{rem}

\end{document}